\documentclass[12pt]{amsart}
\usepackage{amssymb} 
\usepackage{latexsym}
\usepackage[all]{xy}
\usepackage{comment,url}

\newcommand{\aaa}{\mathfrak{a}}
\newcommand{\bbb}{\mathfrak{b}}
\newcommand{\Br}{\operatorname{Br}}
\newcommand{\eps}{{\varepsilon}}
\newcommand{\F}{{\mathbb F}}

\newcommand{\divv}{\operatorname{div}}

\newcommand{\im}{\operatorname{Im}}
\newcommand{\Gm}{{\mathbb G}_m}
\newcommand{\Kbar}{\overline{K}}
\newcommand{\isom}{ \cong }
\newcommand{\tors}{\operatorname{tors}}
\newcommand{\phihat}{\widehat{\phi}}
\newcommand{\psihat}{\widehat{\psi}}
\newcommand{\Div}{\operatorname{Div}}

\newcommand{\Gal}{\operatorname{Gal}}
\newcommand{\Pic}{\operatorname{Pic}}
\newcommand{\PP}{{\mathbb P}}
\newcommand{\Q}{{\mathbb Q}}
\newcommand{\Res}{{\operatorname{Res}}}
\newcommand{\rank}{\operatorname{rank}}
\newcommand{\ra}{{\longrightarrow}}
\newcommand{\Sel}{\operatorname{Sel}}
\newcommand{\Z}{{\mathbb Z}}
\newcommand{\Sha}{\mbox{\wncyr Sh}}

\newfont{\wncyr}{wncyr10 at 12pt}
\newfont{\wncyrten}{wncyr10 at 10pt}

\newenvironment{Proof}{\par\noindent{\sc Proof:}}%
                      {\hspace*{\fill}\nobreak$\Box$\par\medskip}
\newenvironment{ProofOf}[1]{\par\noindent{\sc Proof of #1:}}%
                       {\hspace*{\fill}\nobreak$\Box$\par\medskip}

\newtheorem{Proposition}{Proposition}[section]
\newtheorem{Theorem}[Proposition]{Theorem}
\newtheorem{Lemma}[Proposition]{Lemma}

\theoremstyle{definition}
\newtheorem{Algorithm}[Proposition]{Algorithm}

\newtheorem{Remark}[Proposition]{Remark}
\newtheorem{Example}[Proposition]{Example}

\begin{document}

\title[Higher descents on an elliptic curve]%
{Higher descents on an elliptic curve \\ with a rational 2-torsion point}

\author{Tom~Fisher}
\address{University of Cambridge,
         DPMMS, Centre for Mathematical Sciences,
         Wilberforce Road, Cambridge CB3 0WB, UK}
\email{T.A.Fisher@dpmms.cam.ac.uk}

\date{10th September 2015}

\begin{abstract}
  Let $E$ be an elliptic curve over a number field $K$. Descent
  calculations on $E$ can be used to find upper bounds for the rank of
  the Mordell-Weil group, and to compute covering curves that assist
  in the search for generators of this group. The general method of
  $4$-descent, developed in the PhD theses of Siksek, Womack and
  Stamminger, has been implemented in Magma (when $K=\Q$) and works
  well for elliptic curves with sufficiently small discriminant. By
  extending work of Bremner and Cassels, we describe the improvements
  that can be made when $E$ has a rational $2$-torsion point. In
  particular, when $E$ has full rational $2$-torsion, we describe a
  method for $8$-descent that is practical for elliptic curves $E/\Q$
  with large discriminant.
\end{abstract}

\maketitle

\renewcommand{\baselinestretch}{1.1}
\renewcommand{\arraystretch}{1.3}
\renewcommand{\theenumi}{\roman{enumi}}

\section{Introduction}
\label{introduction}

Let $E$ be an elliptic curve over a number field $K$. For each integer
$n \ge 2$ there is a short exact sequence of abelian groups
\begin{equation*}
  0 \to E(K)/ n E(K) \to S^{(n)}(E/K) \to \Sha(E/K)[n] \to 0.
\end{equation*}
The $n$-Selmer group $S^{(n)}(E/K)$ is finite and effectively
computable.  It gives information about both the Mordell-Weil group
$E(K)$ and the Tate-Shafarevich group $\Sha(E/K)$. The elements of
$S^{(n)}(E/K)$ may be interpreted geometrically as $n$-coverings of
$E$. An {\em $n$-covering} of $E$ is a pair $(C,\pi)$, where $C/K$ is
a smooth curve of genus one, and $\pi : C \to E$ is a morphism defined
over $K$, that fits in a commutative diagram
\[ \xymatrix{ C \ar[d]_{\isom} \ar[dr]^\pi \\ E \ar[r]_{[n]} & E } \]
where the vertical map is an isomorphism defined over $\Kbar$.  The
Selmer group $S^{(n)}(E/K)$ consists of those $n$-coverings $(C,\pi)$
that are everywhere locally soluble, i.e. $C(K_v) \not= \emptyset$ for
all places $v$ of $K$. The subset of those $n$-coverings with $C(K)
\not= \emptyset$ form the image of $E(K)/nE(K)$.  Thus if $C/K$ is a
counter-example to the Hasse Principle then it represents a
non-trivial element of $\Sha(E/K)$.

A {\em first descent} computes the group $S^{(n)}(E/K)$ and represents
its elements as pairs $(C,\pi)$. To compute the group $E(K)/nE(K)$,
and hence the rank of $E(K)$, we must decide which of these
$n$-coverings has a rational point, that is, a point with co-ordinates
in $K$.  Unfortunately there is no known algorithm guaranteed to
determine whether a genus one curve $C/K$ has a rational point. In
practice one starts by searching for rational points of small
height. If no points are found, this might be because $C$ has no
rational points, or because the rational points all have large height.
We attempt to distinguish these two cases by a second descent.

Taking Galois cohomology of the short exact sequence
\[ 0 \ra E[n] \ra E[n^2] \ra E[n] \ra 0, \]
and restricting to Selmer groups, gives an exact sequence
\[ E(K)[n] \ra S^{(n)}(E/K) \ra S^{(n^2)}(E/K) 
\stackrel{\alpha}{\ra} S^{(n)}(E/K). \]
There are then inclusions
\[ E(K)/nE(K) \subset \im(\alpha) \subset  S^{(n)}(E/K). \]
Moreover, the image of  $\alpha$ is the kernel of the 
Cassels-Tate pairing
\[ S^{(n)}(E/K) \times S^{(n)}(E/K) \to \Q/\Z. \] If $\alpha$ maps
$(C_2,\nu_2)$ to $(C_1,\pi_1)$ then $\nu_2$ factors via $\pi_1$ to
give a commutative diagram
\[ \xymatrix{ C_2 \ar[r]^{\pi_2} \ar[d]_{\isom} 
& C_1 \ar[d]_{\isom} \ar[dr]^{\pi_1} \\ 
E \ar[r]_{[n]} & E \ar[r]_{[n]} & E } \]
where the vertical maps are isomorphisms defined over $\Kbar$.

A {\em second descent} computes the fibre of $\alpha$ above
$(C_1,\pi_1)$ and represents its elements as pairs $(C_2, \pi_2)$. If
the fibre is empty then $C_1(K) = \emptyset$.  Otherwise the fibre is
a coset of the image of $S^{(n)}(E/K)$.  We can then try searching for
rational points on each of the genus one curves $C_2$. If we still do
not find a rational point then a third descent may be attempted, and
so on.

More generally, if $\phi : E \to E'$ is an isogeny of degree $n$,
and $\phihat : E' \to E$ is the dual isogeny, then there are
exact sequences
\begin{equation}
\label{kummer}
0 \to E'(K)/ \phi E(K) \to S^{(\phi)}(E/K) \to \Sha(E/K)[\phi] \to 0
\end{equation}
and 
\[ E(K)[\phi] \ra S^{(\phihat)}(E'/K) \ra S^{(n)}(E'/K) 
\stackrel{\alpha}{\ra} S^{(\phi)}(E/K). \]
Moreover, the image of $\alpha$ is the kernel of the 
Cassels-Tate pairing 
\[ S^{(\phi)}(E/K) \times S^{(\phi)}(E/K) \to \Q/\Z. \]
The terminology of first and second descents carries over as before.

In this paper we are concerned with $\phi : E \to E'$ an isogeny of
degree $2$. Thus our work applies to any elliptic curve with a
rational $2$-torsion point.  The first descent in this case is descent
by $2$-isogeny.  This is very well known; see for example \cite[\S
14]{CaL} or \cite[X.4.9]{Si1}. The second descent is described in
\cite[\S 5]{BSDII} and \cite{Cr}, although in neither case do
the authors claim any particular originality. This corresponds to a
$2$-descent on $E$. The starting point for our work is the paper of
Bremner and Cassels~\cite{BC} that carries out the third and fourth
descents for elliptic curves $E/\Q$ of the form $y^2=x(x^2-4p)$, where
$p$ is a prime with $p \equiv 5 \pmod{8}$. This corresponds to a
$4$-descent on $E$.  As observed by Siksek~\cite[\S 4.6]{Siksek},
their method can be applied more generally. 

The general method of $4$-descent, developed in \cite{Siksek},
\cite{MSS}, \cite{Womack}, \cite{Stamminger}, (see also
\cite{improve4}), requires that we compute the class group and units
of a degree $4$ extension of $K$.  This has been implemented in Magma
\cite{magma} (when $K=\Q$) and works well for elliptic curves with
sufficiently small discriminant. In contrast the method of Bremner and
Cassels (specific to elliptic curves with a rational $2$-torsion
point) requires no class group and unit calculations, beyond those for
the field $K$ itself.  Instead the global part of the calculation
requires that we solve conics over $K$, and over quadratic extensions
of $K$.

We make the following improvements.
\begin{itemize}
\item We use the Cassels-Tate pairing to efficiently compute 
upper bounds for the rank. This was not required in \cite{BC}, 
since for the curves considered there, descent by $2$-isogeny 
already shows that the rank is at most $1$.
\item We extend to a fifth descent, and in cases where $E$ 
has full rational $2$-torsion, also a sixth descent. The
latter corresponds to $8$-descent on $E$.
\item We replace the problem of solving a conic over a quadratic
extension of $K$, with that of solving a quadric surface over $K$.
Taking $K= \Q$ the latter can be solved efficiently using 
an algorithm of Simon \cite{SimonQuadrics}.
\end{itemize}

As discussed further in Section~\ref{sec:examples}, the motivation for
our work came from the desire to find curves of large rank in certain
families of elliptic curves over $\Q$, for example the family of
elliptic curves with a given torsion subgroup (of even order), or the
quadratic twists of a given elliptic curve with a rational $2$-torsion
point. The methods we describe can be used to quickly eliminate many
curves which on the basis of (say) descent by 2-isogeny appear to be
candidates for large rank, but which instead have large $2$-primary
part of $\Sha$. 

The paper is organised as follows. In the first three sections we let
$\phi: E \to E'$ be any isogeny of prime degree $p$.  In
Section~\ref{intro} we introduce the higher descent pairings we use to
bound the rank of an elliptic curve. Then in Sections~\ref{sec:ctp}
and~\ref{sec:book} we explain how these pairings are related to the
Cassels-Tate pairing, and outline our methods for computing them. In
Section~\ref{sec:outline} we give a short self-contained account of
$4$-descent on an elliptic curve with a rational $2$-torsion
point. This is then related to the work of Bremner and Cassels in
Section~\ref{sec:compare}. In the next two sections we describe our
refinements for carrying out the fifth and sixth descents. In
Section~\ref{sec:res} we explain how to replace the conics over
quadratic extensions by quadric surfaces. Finally in
Section~\ref{sec:examples} we give some examples.

Our implementation of the methods described in this paper (in the case
$K= \Q$) has been contributed to Magma (version 2.21) and
is available via the function {\tt TwoPowerIsogenyDescentRankBound}.

\section{Higher descent pairings}
\label{intro}

Let $\phi : E \to E'$ be an isogeny of elliptic curves 
defined over a number field $K$. We suppose that $\deg \phi = p$
is a prime. By \cite[X.4.7]{Si1} there
is an exact sequence
\begin{align*} 
0 \ra E(K)&[\phi] \ra E(K)[p] \ra E'(K)[\phihat] \ra \\
  &\ra E'(K)/\phi E(K) \ra E(K)/p E(K) \ra E(K)/\phihat E'(K) \ra 0.
\end{align*}
Writing $\dim$ for the dimension of an $\F_p$-vector 
space, it follows that
\[ \rank E(K) 
= \dim \frac{E'(K)}{\phi E(K)}+ \dim \frac{E(K)}{\phihat E'(K)}
  - \dim E(K)[\phi] - \dim E'(K)[\phihat]. \]

Let $S_1 = S^{(\phi)}(E/K)$ and  $S_1' = S^{(\phihat)}(E'/K)$ 
be the Selmer groups attached to the isogenies $\phi$ and $\phihat$.
More generally let $S_m \subset S_1$ be the image of $S^{(p^n)}(E'/K)$ 
if $m = 2n$ is even, and the image of 
$S^{(p^n\phi)}(E/K)$ if $m = 2n+1$ is odd.
The subspaces $S'_m \subset S'_1$ are defined in the same way,
after swapping the roles of $E$ and~$E'$.
There are inclusions of $\F_p$-vector spaces
\begin{equation}
\label{filt1}
 \frac{E'(K)}{\phi E(K)} \subset \ldots \subset S_3 \subset S_2
\subset S_1 = S^{(\phi)}(E/K) 
\end{equation}
and 
\begin{equation}
\label{filt2}
 \frac{E(K)}{\phihat E'(K)} \subset \ldots \subset S'_3 \subset S'_2
\subset S'_1 = S^{(\phihat)}(E'/K).
\end{equation}

As we prove in Section~\ref{sec:ctp}, the Cassels-Tate pairing 
induces the following pairings of $\F_p$-vector spaces.

\begin{Theorem}
\label{thm:pair}
Let $m \ge 1$ be an integer.
\begin{enumerate}
\item
 If $m$ is odd then there are alternating pairings
\[ \Theta_m : S_m \times S_m \to \F_p \qquad \text{ and } \qquad 
\Theta'_m : S'_m \times S'_m \to \F_p  \]
with kernels $S_{m+1}$ and $S'_{m+1}$. 
\item
If $m$ is even then there is a pairing 
\[ \Theta_m : S_m \times S'_m \to \F_p \]
with left kernel $S_{m+1}$ and right kernel $S'_{m+1}$. 
\end{enumerate}
\end{Theorem}

It is clear that each time we compute one of the pairings
$\Theta_m$ or $\Theta_m'$ our upper bound for 
the rank of $E(K)$ either stays the same (if the pairing is identically
zero) or decreases by an even integer. The following lemma
is useful for comparing our bounds on the rank with those 
obtained by other methods.

\begin{Lemma} 
\label{lem:pn}
The upper bound on the rank of $E(K)$ obtained by $p^n$-descent is
\begin{equation*}
\rank E(K) \le \dim S_{2n-1} + \dim S'_{2n} - 
 \dim E(K)[\phi] - \dim E'(K)[\phihat].
\end{equation*}
\end{Lemma}
\begin{Proof}
From the commutative diagram with exact rows
\[ \xymatrix{ 
E'(K)[\phihat] \ar[r] \ar@{=}[d] & S^{(p^{n-1}\phi)}(E/K) \ar[r] \ar[d] 
& S^{(p^n)}(E/K) \ar[r] \ar[d]^{\beta} & S^{(\phihat)}(E'/K)  \ar@{=}[d] \\
E'(K)[\phihat] \ar[r] & S^{(\phi)}(E/K) \ar[r] &
S^{(p)}(E/K) \ar[r] & S^{(\phihat)}(E'/K)
} \]
we obtain an exact sequence
\[ 0 \ra E(K) [\phi] \ra E(K)[p] \ra E'(K)[\phihat]
 \ra S_{2n-1} \ra \im(\beta) \ra S'_{2n} \ra 0. \]
Therefore the upper bound obtained by $p^n$-descent,
\[ \rank E(K) \le \dim \im(\beta) - \dim E(K)[p], \]
is the same as that in the statement of the lemma.
\end{Proof}

The filtrations~(\ref{filt1}) and~(\ref{filt2}) also give
information about the Tate-Shafarevich groups of $E$ and $E'$. Let
$\Sha_m \subset \Sha_1 = \Sha(E/K)[\phi]$ be the image of
$\Sha(E'/K)[p^n]$ if $m = 2n$ is even, and the image of
$\Sha(E/K)[p^n\phi]$ if $m = 2n+1$ is odd.  The subspaces $\Sha'_m
\subset \Sha'_1 = \Sha(E'/K)[\phihat]$ are defined in the same way,
after swapping the roles of $E$ and~$E'$. For each integer 
$m \ge 1$ we have short exact sequences 
\[  0 \ra E'(K)/\phi E(K) \ra S_m \ra \Sha_m \ra 0 \]
and 
\[  0 \ra E(K)/\phihat E'(K) \ra S'_m \ra \Sha'_m \ra 0 \rlap{.} \]
In situations where we succeed in computing the rank of $E(K)$ we
have $\Sha_m = \Sha'_m = 0$ for all $m$ sufficiently large.

\begin{Remark} 
\label{rem:getsha}
Suppose we have computed $\Sha_m$ and $\Sha'_m$ 
for all $m \ge 1$. Then from the exact sequences
\begin{align*}
&0 \ra \Sha(E'/K)[p^{n-1} \phihat] \ra \Sha(E'/K)[p^n] \ra \Sha_{2n} 
\ra 0, \\
&\, 0 \ra \Sha(E/K)[p^{n}]  \ra \Sha(E/K)[p^n \phi] \ra \Sha_{2n+1} 
\ra 0,
\end{align*}
and their analogues for $\Sha'_m$, we can read off the orders of
$\Sha(E/K)[p^n]$ and $\Sha(E'/K)[p^n]$ for all $n \ge 1$. 
This information determines the group structure of the 
$p$-primary parts of $\Sha(E/K)$ and $\Sha(E'/K)$. 
\end{Remark}

In this paper we take $p=2$.  We show how to compute $S_m$ and $S'_m$
for $m \le 5$, by a method whose global part only requires that we
solve quadratic forms of ranks $3$ and $4$ over~$K$.  When $E$ has
full rational $2$-torsion we also compute $S'_6$.  By
Lemma~\ref{lem:pn} the upper bound $\rank E(K) \le \dim S_5 + \dim
S'_6 - 2$ is then the same as that obtained by $8$-descent on $E$.

\section{The Cassels-Tate pairing}
\label{sec:ctp}

Let $E$ be an elliptic curve over a number field $K$.
The Cassels-Tate pairing is an alternating bilinear pairing
\[\langle~,~\rangle : \Sha(E/K) \times \Sha(E/K) \to \Q/\Z.\]
It has the following properties.

\begin{Theorem}
\label{thm:sha}
Let $\psi : C \to D$ be an isogeny of elliptic curves over $K$.
\begin{enumerate}
\item $\langle \psi x,y \rangle  = \langle x, \psihat y \rangle$
for all $x \in \Sha(C/K)$ and $y \in \Sha(D/K)$.
\item $x \in \Sha(D/K)$ belongs to the image of $\psi : \Sha(C/K) 
\to \Sha(D/K)$ if and only if $\langle x,y \rangle = 0$ for all
$y$ in the kernel of $\psihat: \Sha(D/K) \to \Sha(C/K)$.
\end{enumerate}
\end{Theorem}
\begin{Proof}
For (i) see \cite[Section 2]{CaVIII}, and for (ii) see 
\cite[Theorem 3]{CTPPS}. These results do {\em not} depend
on finiteness of $\Sha$.
\end{Proof}

\begin{ProofOf}{Theorem~\ref{thm:pair}}
We keep the notation of Section~\ref{intro}, except that
the pairings $\Theta_m$ and $\Theta_m'$ will take values in 
$\frac{1}{p}\Z/\Z$ instead of $\F_p$. We make frequent implicit
use of the exact sequence~(\ref{kummer}). In particular it makes
sense to evaluate the Cassels-Tate pairing on Selmer group
elements. \\
(i) Let $m = 2n + 1$. Let $\xi,\eta \in S_m$, and let 
$\xi_1, \eta_1 \in S^{(p^n \phi)}(E/K)$ with 
$\xi_1 \mapsto \xi$ and $\eta_1 \mapsto \eta$. We define
$\Theta_m(\xi,\eta) = \langle \xi_1, \eta \rangle
 = \langle \xi, \eta_1 \rangle$.
These last two expressions are equal by Theorem~\ref{thm:sha}(i)
with $\psi = p^n$. Therefore $\Theta_m(\xi,\eta)$
is independent of the choices of $\xi_1$ and $\eta_1$.
By Theorem~\ref{thm:sha}(ii) the kernel of $\Theta_m$ is
\[  \{ \xi \in S_m \, | \, \langle \xi,\eta_1 \rangle = 0 
\text{ for all } \eta_1 \in S^{(p^n \phi)}(E/K) \} = S_{m+1}. \]
Now let $\psi = p^{n/2}$ or $p^{(n-1)/2} \phi$ according as $n$ is even
or odd. By Theorem~\ref{thm:sha}(i) we have 
\[ \Theta_m(\xi,\xi) = \langle \xi_1, \psihat \psi \xi_1 \rangle 
= \langle \psi \xi_1, \psi \xi_1 \rangle = 0. \]
Therefore $\Theta_m$ is alternating. The definition and properties of 
$\Theta'_m$ are obtained in the same way, after swapping the roles 
of $E$ and $E'$. \\
(ii) Let $m = 2n$. Let $\xi \in S_m$ and $\eta \in S'_m$, and let 
$\xi_1 \in S^{(p^n)}(E'/K)$ and $\eta_1 \in S^{(p^n)}(E/K)$ 
with $\xi_1 \mapsto \xi$ and $\eta_1 \mapsto \eta$. We define
$\Theta_m(\xi,\eta) = \langle \xi_1, \eta \rangle 
 = \langle \xi, \eta_1 \rangle$.
These last two expressions are equal by 
Theorem~\ref{thm:sha}(i) with $\psi = p^{n-1} \phi$. 
Therefore $\Theta_m(\xi,\eta)$ is independent of the 
choices of $\xi_1$ and $\eta_1$.
By Theorem~\ref{thm:sha}(ii) the left kernel of $\Theta_m$ is
\[ 
\{ \xi \in S_m \, | \, \langle \xi,\eta_1 \rangle = 0 
\text{ for all } \eta_1 \in S^{(p^n)}(E/K) \} = S_{m+1}. \]
Likewise the right kernel is $S'_{m+1}$. 
\end{ProofOf}

Next we describe a method for computing the Cassels-Tate pairing.
Suppose that $E[\phi] \isom \mu_p$ and $E'[\phihat] \isom
\Z/p\Z$ as Galois modules. Then  
$H^1(K,E[\phi]) = K^\times/(K^\times)^p$ and 
$H^2(K,E[\phi]) = \Br(K)[p]$. If $\psi : E' \to F$ is an isogeny
defined over~$K$ then from the short exact sequence
\[ 0 \to E[\phi] \to E[\psi \phi] \to E'[\psi] \to 0 \]
we obtain a long exact sequence
\begin{equation}
\label{les}
 \ldots \to K^\times/(K^\times)^p \to H^1(K,E[\psi \phi]) \to
H^1(K,E'[\psi]) \to  \Br(K)[p] \to \ldots 
\end{equation}

Let $x \in S^{(\psi)}(E'/K)$ and $y \in S^{(\phihat)}(E'/K)$. By the 
local-to-global principle for the Brauer group, $x$ lifts to 
$x' \in H^1(K,E[\psi \phi])$. Let $C$ and $D_1$ be the covering
curves corresponding to $x$ and $x'$. These fit in a commutative
diagram
\begin{equation*}
\xymatrix{ D_1 \ar[d]_\isom \ar[r]^{\pi} & C \ar[d]_\isom \ar[dr] \\
E \ar[r]^\phi & E' \ar[r]^\psi & F }
\end{equation*}
where the vertical maps are isomorphisms defined over $\Kbar$,
but all other maps are morphisms defined over $K$.

Let $T$ be a generator for $E'(K)[\phihat] \isom \Z/p\Z$, and let
$\bbb \in \Div^0 (C)$ correspond to $T$ under the isomorphism 
of Galois modules
$\Pic^0(C) \isom E'$. Since $T \in E'(K)$ and $C$ is everywhere
locally soluble, we may choose $\bbb$ to be $K$-rational.
Then there exists $f \in K(C)$ with $\divv(f) = p \bbb$.
We say that $f$ is a {\em pushout function}.
If we scale $f$ suitably then $\pi^* f = g^p$ for some 
$g \in K(D_1)$. In other words,
if we identify $K(C)$ as a subfield of 
$K(D_1)$ via pull-back by $\pi$, then $K(D_1) = K(C)(\sqrt[p]{f})$.

For each place $v$ of $K$ there is a local pairing 
\begin{equation}
\label{locpair}
 (~,~)_v : H^1(K_v,E[\phi]) \times H^1(K_v,E'[\phihat]) 
\to \tfrac{1}{p} \Z/\Z 
\end{equation}
given by the Weil pairing, cup product and the local invariant map.
We identify $H^1(K_v,E[\phi]) = K_v^\times/(K_v^\times)^p$.
The Cassels-Tate pairing is given by
\[ \langle x, y \rangle = \sum_v ( f(P_v),y)_v \]
where for each place $v$ of $K$ we choose a local point $P_v \in C(K_v)$ 
avoiding the zeros and poles of $f$.
This is a sum over all places of $K$, but in fact there is no
contribution from the primes outside a finite set of primes
where $C$ and $f$ have bad reduction.

The Selmer group attached to $\phi$ is
\[ S^{(\phi)}(E/K) = \{ \xi \in K^\times/(K^\times)^p \, | \, 
\xi \in \im(\delta_{\phi,v}) \text{ for all places } v \, \} \]
where $\delta_{\phi,v} : E'(K_v) \to K_v^\times/(K_v^\times)^p$ is the local
connecting map. 

Let $D_{\xi}$ be the covering of $C$ that is $K$-birational
to $\{ f(P) = \xi z^p \} \subset C \times \Gm$.  
The Selmer set associated to the pair $(C,f)$ is 
\begin{align*}
  S(C,f) &= \{ \, \xi \in K^\times/(K^\times)^p \, | \, 
 D_{\xi}(K_v) \not= \emptyset \text{ for all places } v \, \} \\
&= \{ \, \xi \in K^\times/(K^\times)^p \, | \, 
\xi \equiv f(P_v) \!\!\! \mod{\im (\delta_{\phi,v})} 
   \text{ for all places } v \, \}.
\end{align*}
The Selmer set tells us which twists of $D_1$ (as a
covering of $C$) are everywhere locally soluble. 
If $\langle x,y\rangle \not= 0$ for some $y \in S^{(\phihat)}(E'/K)$ 
then there are no such twists, and the Selmer set $S(C,f)$ is empty. 
Otherwise  $S(C,f)$ is a coset of $S^{(\phi)}(E/K)$ 
in $K^\times/(K^\times)^p$. If we have already computed the pairing
(using the formula above) then we can solve for a coset 
representative by linear algebra over $\F_p$. 

In this paper we take $p=2$. Then $E[\phi] \isom E'[\phihat] \isom
\Z/2\Z$ and the local pairing~(\ref{locpair}) is the usual
(i.e. quadratic) Hilbert norm residue symbol. As above, let $C$ be the
covering curve corresponding to $x \in S^{(\psi)}(E'/K)$. Suppose we
know equations for~$C$. Then by computing a pushout function $f$ on
$C$, and using it to compute the Cassels-Tate pairing, we can
determine whether $x$ lifts to $x' \in S^{(\psi \phi)}(E/K)$. If it
does lift then the covering curve corresponding to $x'$ is of the form
$D_\xi$ for some $\xi \in S(C,f)$.  In particular we know equations
for $D_\xi$.  We can then replace $\psi$ by $\psi \phi$, swap the
roles of $E$ and $E'$, and repeat.

At each stage we have a basis for $S_m \subset K^\times/(K^\times)^2$,
and for each basis element $\xi$, equations for an everywhere locally
soluble $2^m$-isogeny covering of $E'$ that factors via the
$\phihat$-covering of $E'$ corresponding to $\xi$.  Using this data we
can compute the pairing $\Theta_m$ in Theorem~\ref{thm:pair}, and
hence the subspace $S_{m+1} \subset S_m$. For the next iteration we
need to compute an everywhere locally soluble $2^{m+1}$-isogeny
covering for each basis element of $S_{m+1}$. It often happens that
these curves are degree-$2$ coverings of some of the curves we already
found.  However, in general a rather subtle iterative procedure is
required. This is the subject of Section~\ref{sec:book}.

The two key issues we must address are the following.
\begin{itemize}
\item How do we find ``nice'' equations for the covering curves?
\item How do we compute the pushout functions?
\end{itemize}
These questions are related, in that a good answer to the first
question helps with answering the second.  The proof that pushout
functions exist uses the local-to-global principle for $\Br(K)[2]$. So
one might expect that the second question comes down to solving conics
over $K$. Indeed if $K$ is a number field, then every element of
$\Br(K)[2]$ can be represented by a conic.  However this last
statement is not true over arbitrary fields, and the proof over number
fields itself uses the local-to-global principle for the Brauer
group. So the best we can say for arbitrary $m$ (see \cite{SD}) is
that the second question reduces to that of finding rational points on
certain Brauer-Severi varieties.

\section{An iterative procedure}
\label{sec:book}

In this section we describe the structure of our program for computing
the subspaces $S_m$ and $S'_m$. The bookkeeping is somewhat involved,
but is needed to get around the fact that the formula we gave in
Section~\ref{sec:ctp}, for computing the Cassels-Tate pairing, only
applies when one of the arguments is killed by the $p$-isogeny $\phi$
or its dual $\phihat$.

First we lighten our notation by writing
\[ \Sel(m) = \left\{ \begin{array}{ll} 
 S^{(p^n)}(E/K) \text{ or } S^{(p^n)}(E'/K) & \text{ if } m = 2n, \\
 S^{(p^n \phi)}(E/K) \text{ or } S^{(p^n \phihat)}(E'/K) & \text{ if } 
m = 2n+1.  \end{array} \right. \]
Which of the two groups we mean is often determined by the context.
For example when we write
\[ \ldots \ra \Sel(3) \ra \Sel(2) \ra \Sel(1) \]
this could either mean 
\[ \ldots\ra S^{(p \phi)}(E/K) \ra S^{(p)}(E'/K) \ra S^{(\phi)}(E/K) \]
or 
\[\ldots \ra S^{(p \phihat)}(E'/K) \ra S^{(p)}(E/K) \ra S^{(\phihat)}(E'/K), \]
whereas
\[  \Sel(1) \ra \Sel(2) \ra \Sel(3) \ra\ldots \]
could either mean
\[ S^{(\phi)}(E/K) \ra S^{(p)}(E/K) \ra S^{(p\phi)}(E/K) \ra  \ldots\]
or 
\[ S^{(\phihat)}(E'/K) \ra S^{(p)}(E'/K) \ra S^{(p\phihat)}(E'/K) \ra \ldots. \]

\medskip

Let $S_m$ be the image of $\Sel(m) \to \Sel(1)$. This $\F_p$-vector
space was denoted $S_m$ or $S'_m$ in Section~\ref{intro}. We also
write $\Theta_m : S_m \times S_m \to \F_p$ for the pairings in
Theorem~\ref{thm:pair}, even though some of the $S_m$ and $\Theta_m$
should really be $S'_m$ and $\Theta'_m$. These abuses of notation 
are made to simplify the results in this section, and (with
the exception of Lemma~\ref{lem:book8}) will not be 
used elsewhere in the paper.

Let $\alpha_1, \beta_1 \in S_m$. Then $\Theta_m(\alpha_1,\beta_1) = 
\langle \alpha_ m, \beta_1 \rangle$ where $\alpha_m \in \Sel(m)$ 
is any lift of $\alpha_1$. In order to compute $\Theta_m$, and hence 
by Theorem~\ref{thm:pair} the subspace $S_{m+1}$, we must describe 
how to lift $\alpha_1$ to $\alpha_m$. In the case $m=3$ it is 
convenient to break this down into the following steps.
\begin{enumerate}
\item Solve for $\alpha_2 \in \Sel(2)$ with $\alpha_2 \mapsto \alpha_1$.
\item Solve for $\xi \in S_1$ with $\langle \xi + 
   \alpha_2 ,\eta \rangle = 0$ for all $\eta \in S_1$. 
\item Solve for $\alpha_3 \in \Sel(3)$ with $\alpha_3 \mapsto \xi + \alpha_2$.
\end{enumerate}
For general $m \ge 1$ we use the following algorithm.

\renewcommand{\theenumi}{\arabic{enumi}}

\begin{Algorithm} 
\label{alg1}
INPUT: $\alpha_1 \in S_m$ and the pairings $\Theta_\ell$ for $\ell < m$. \\
OUTPUT: $\alpha_m \in \Sel(m)$ with $\alpha_m \mapsto \alpha_1$.
\begin{enumerate}
\item {\tt if $m$ = $1$ then return $\alpha_1$; end if}
\item {\tt $k \leftarrow 1$}
\item {\tt $t_1 \leftarrow m - 1$} 
\item {\tt while true do}
\item ~\indent {\tt solve for $\alpha_{k+1} \in \Sel(k+1)$ with $\alpha_{k+1}
\mapsto \alpha_k$}
\item ~\indent $k \leftarrow k + 1$
\item ~\indent {\tt if $k=m$ then return $\alpha_m$; end if}
\item ~\indent $t_k \leftarrow 1$
\item ~\indent {\tt while $t_{k-1} = t_k$ do}
\item ~\indent\indent $k \leftarrow k - 1$
\item ~\indent\indent $t_k \leftarrow t_k + 1$
\item ~\indent {\tt end while}
\item ~\indent $\ell \leftarrow t_k$
\item ~\indent {\tt solve for $\xi \in S_\ell$ with $\Theta_\ell(\xi,\eta)
 + \langle \alpha_{k+\ell-1} , \eta \rangle = 0$ for all $\eta \in S_\ell$.} 
\item ~\indent $\alpha_k \leftarrow \xi + \alpha_k$
\item {\tt end while}
\end{enumerate}
\end{Algorithm}
 
\renewcommand{\theenumi}{\roman{enumi}}

The next two lemmas are used to show that Algorithm~\ref{alg1} is
correct.

\begin{Lemma}
\label{lem:book1}
Suppose $\alpha_k \in \Sel(k)$ lifts to $\alpha_{k+\ell-1} \in \Sel(k+\ell-1)$.
Let $\xi \in \Sel(1)$. The following are equivalent.
\begin{enumerate}
\item $\langle \xi + \alpha_k, \beta \rangle = 0$ for all 
$\beta \in \Sel(\ell)$, 
\item $\xi \in S_\ell$ and $\Theta_\ell(\xi,\eta)
 + \langle \alpha_{k+\ell-1} , \eta \rangle = 0$ for all $\eta \in S_\ell$. 
\end{enumerate}
Moreover, if $\alpha_k \mapsto \alpha_{k-1} \in \Sel(k-1)$ and 
$\langle \alpha_{k-1}, \beta \rangle =0$ for all $\beta \in \Sel(\ell+1)$,
then (i) and (ii) hold for some $\xi \in \Sel(1)$.
\end{Lemma}

\begin{Proof}
  We make repeated use of Theorem~\ref{thm:sha}.  Since $\alpha_k$
  lifts to $\Sel(k+ \ell - 1)$ we have $\langle \alpha_k, \beta
  \rangle = 0$ for all $\beta \in \Sel(\ell - 1)$.  Assuming (i) we
  have $\langle \xi, \beta \rangle = 0$ for all $\beta \in \Sel(\ell -
  1)$, and therefore $\xi \in S_\ell$.  If $\beta \in \Sel(\ell)$ and
  $\eta \in S_\ell$ with $\beta \mapsto \eta$ then
\[  \langle \xi + \alpha_k , \beta \rangle 
=  \langle \xi, \beta \rangle  + \langle \alpha_k , \beta \rangle 
= \Theta_\ell(\xi,\eta) +  \langle \alpha_{k+\ell-1} , \eta \rangle.  
\]
This proves the equivalence of (i) and (ii). For the last part,
Theorem~\ref{thm:sha}(ii) shows there exists $\alpha_{k + \ell} 
\in \Sel(k+\ell)$ with $\alpha_{k+\ell} \mapsto \alpha_{k-1}$. 
Let $\widetilde{\alpha}_k \in \Sel(k)$ be the image of 
$\alpha_{k+\ell}$. By Theorem~\ref{thm:sha}(i) we have  
$\langle \widetilde{\alpha}_k , \beta \rangle =0$ for all $\beta
\in \Sel(\ell)$. Since $\widetilde{\alpha}_k$ and $\alpha_k$ have
the same image in $\Sel(k-1)$ we have $\widetilde{\alpha}_k = 
\xi + \alpha_k$ for some $\xi \in \Sel(1)$. This proves (i).
\end{Proof}

\begin{Lemma}
\label{lem:book2}
At the start and end of the main loop in Algorithm~\ref{alg1},
we have $\alpha_k \mapsto \alpha_{k-1} \mapsto \ldots \mapsto 
\alpha_2 \mapsto \alpha_1$ and
\begin{equation}
\label{tprop}
\text{$\langle \alpha_j, \beta \rangle = 0$ for all $\beta 
\in \Sel(t_j)$ and $1 \le j \le k$.}
\end{equation}
\end{Lemma} 
\begin{Proof}
Initially we have $k=1$ and $t_1 = m-1$. Since $\alpha_1$ is in 
$S_m$ it satisfies $\langle \alpha_1, \beta \rangle = 0$ for all
$\beta \in \Sel(m-1)$. At the end of each loop it suffices to 
prove~(\ref{tprop}) with $j=k$, since the cases $j<k$ carry over
from the previous iteration. It is easy to check that 
each time we reach line 14, we have $k \ge 2$ and $t_{k-1} > t_k = \ell$.
Moreover $\alpha_{k+ \ell -1}$ is the element computed in line 5.
By~(\ref{tprop}) from the previous iteration we have 
$\langle \alpha_{k-1},\beta \rangle = 0$ for all $\beta \in \Sel(\ell+1)$.
Lemma~\ref{lem:book1} shows that there exists $\xi \in S_\ell$ satisfying
the condition on line 14, and that after modifying $\alpha_k$ on 
line 15, the new $\alpha_k$ satisfies~(\ref{tprop}) with $j=k$.
\end{Proof}

As a special case of Lemma~\ref{lem:book2} we have $\langle
\alpha_k, \beta \rangle =0$ for all $\beta \in \Sel(t_k)$. Since
$t_k \ge 1$, the element
$\alpha_{k+1}$ in line 5 exists by Theorem~\ref{thm:sha}(ii).
This completes the proof that Algorithm~\ref{alg1} is correct.

The reader is encouraged to write out the sequence $t_1, \ldots, t_k$
at each stage of the algorithm (for some small values of $m$). It is
not hard to show that if $m \ge 2$ then we start the main loop exactly
$2^{m-2}$ times. This exponential growth does not concern us much,
since (even in the case $p=2$) the computation of pushout functions
(required in lines 5 and 14, as explained in Section~\ref{sec:ctp}) is
only practical for small values of $m$.

In practice we can often take $\xi = 0$ in line 14. In such cases we
may already know $\alpha_{k+1} \in \Sel(k+1)$ with $\alpha_{k+1}
\mapsto \alpha_k$. There is then no need to recompute $\alpha_{k+1}$
in line 5.

\section{Computing $4$-coverings}
\label{sec:outline}

In this section we give a brief self-contained account of $4$-descent
on an elliptic curve with a rational $2$-torsion point. As we discuss
in Section~\ref{sec:compare}, it is based on work of Bremner and
Cassels \cite{BC}. See the introduction for further references.

Let $E$ be an elliptic curve over a field $K$ with
a rational 2-torsion point $T$, say 
\begin{equation}
\label{weqnE}
E: \quad y^2 = x(x^2 + ax + b) \qquad \quad  T = (0,0) 
\end{equation}
for some $a,b \in K$. The discriminant condition is $2 b (a^2 - 4 b)
\not= 0$. Let $\phi : E \to E'$ be the $2$-isogeny with kernel
$\{0,T\}$ and let $\phihat : E' \to E$ be the dual isogeny. The
connecting map $\delta : E(K)/\phihat E'(K) \to K^\times /
(K^\times)^2$ is given by $P = (x,y) \mapsto x$ for all $ P \not= 0,
T$.

Suppose $P=(x,y) \in E(K)$ with $\delta(P) = \xi_1 \mod{(K^\times)^2}$. 
Then $x = \xi_1 (s/t)^2$ and $y = \xi_1(rs/t^3)$ where
\begin{equation}
\label{cov1}
 r^2 = \xi_1 s^4 + a s^2 t^2 + (b/\xi_1) t^4. 
\end{equation}
Parametrising a conic over $K$ gives 
\begin{equation}
\label{param1}
(s^2 : t^2 : r) = (f(l,m) :g(l,m) : h(l,m)) 
\end{equation}
where $f$, $g$ and $h$ are binary quadratic forms. Then 
\begin{equation}
\label{cov2}
 f(l,m) = \xi_2 s^2 \quad \text{ and } \quad g(l,m) = \xi_2 t^2 
\end{equation}
for some $\xi_2 \in K^\times$. Parametrising each of these conics 
over $K$ gives
\begin{align*}
( l : m : s ) & = ( p_1(c,d) : p_2(c,d) : p_3(c,d) ) \\
( l : m : t ) & = ( q_1(\theta,\psi) : q_2(\theta,\psi) : q_3(\theta,\psi) ) 
\end{align*}
where the $p_i$ and $q_i$ are binary quadratic forms.
Then 
\begin{equation}
\label{cov3}
   p_1(c,d) = \xi_3 q_1(\theta,\psi) \quad \text{ and } \quad 
   p_2(c,d) = \xi_3 q_2(\theta,\psi) 
\end{equation}
for some $\xi_3 \in K^\times$. 

Let $L = K(\sqrt{b/\xi_1})$ and write $f(l,m) = 
\kappa (l - \eps m)(l - \overline{\eps}m)$ where $\kappa \in K$ 
and $\eps, \overline{\eps} \in L$. Then 
\begin{equation*}
 \kappa (p_1 - \eps p_2) (p_1 - \overline{\eps}p_2) = \xi_2 p_3^2. 
\end{equation*}
Since $p_1$ and $p_2$ are coprime in $K[c,d]$ it follows that
\begin{equation*}
 p_1(c,d) - \eps p_2(c,d) = \xi_3 \alpha ( c + \gamma d )^2 
\end{equation*}
for some $\alpha,\gamma \in L$. Hence by~(\ref{cov3}) we have 
\begin{equation*}
 q_1(\theta,\psi) - \eps q_2(\theta,\psi) = \alpha ( c + \gamma d )^2. 
\end{equation*}
Parametrising a conic over $L$ gives
\[ (\theta : \psi : c + \gamma d ) 
  = (Q_1(\lambda,\mu) : Q_2(\lambda,\mu)  : Q_3(\lambda,\mu) ) \]
where $Q_1$, $Q_2$ and $Q_3$ are binary quadratic forms. Then 
\begin{equation}
\label{eqnL}
\theta = \pi Q_1(\lambda,\mu) \quad \text{ and } \quad 
   \psi = \pi Q_2(\lambda,\mu) 
\end{equation}
for some $\pi \in L^\times$. Let $1, \beta$ be a basis for $L$ over $K$.
Writing $\lambda = x + \beta y$ and $\mu = u + \beta v$ we expand
to give
\begin{align*}
\pi Q_1(\lambda,\mu) &= F_1(x,y,u,v) + \beta F_2(x,y,u,v) \\
\pi Q_2(\lambda,\mu) &= G_1(x,y,u,v) + \beta G_2(x,y,u,v) 
\end{align*}
where $F_1$, $F_2$, $G_1$ and $G_2$ are quadratic forms with coefficients
in $K$. Since $\theta,\psi \in K$ it follows by~(\ref{eqnL}) that
\begin{equation}
\label{finalqi}
F_2(x,y,u,v) = G_2(x,y,u,v) = 0. 
\end{equation}
Up to linear changes of co-ordinates defined over $K$, this quadric
intersection depends only on the image of $\pi$ in 
$L^\times/K^\times (L^\times)^2$, and hence only on 
$\xi_4 := N_{L/K}(\pi) \in K^\times/(K^\times)^2$. We may recover $\pi$
from $\xi_4$ by solving a conic over $K$.

In the case $K$ is a number field, there are at each stage
($j=1,2,3,4$) only finitely many $\xi_j \in K^\times/(K^\times)^2$
for which the corresponding covering curves are everywhere
locally soluble. Thus each rational point on $E$ lifts 
to one of only finitely many quadric intersections~(\ref{finalqi}).

The equations~(\ref{cov1}), (\ref{cov2}), (\ref{cov3}),
(\ref{finalqi}) define covering curves $C_1, \ldots, C_4$ that fit in
a commutative diagram
\begin{equation}
\label{covdiag}
\begin{aligned}
  \xymatrix{ C_4 \ar[d]_\isom \ar[r]^{\pi_4} & C_3 \ar[d]_\isom
    \ar[r]^{\pi_3}
    & C_2 \ar[d]_\isom \ar[r]^{\pi_2} & C_1 \ar[d]_\isom \ar[dr]^{\pi_1} \\
    E \ar[r]^\phi & E' \ar[r]^\phihat & E \ar[r]^\phi & E'
    \ar[r]^\phihat & E }
\end{aligned}
\end{equation}
where the vertical maps are isomorphisms defined over $\Kbar$, and all
other maps are morphisms of degree $2$ defined over $K$. Explicitly
the $C_j$ have equations:
\[ C_1 = \{ r^2 = \xi_1 s^4 + a s^2 t^2 + (b/\xi_1) t^4 \} = \{ y^2 =
f(l,m) g(l,m) \}, \]
\[ C_2 = \left\{ 
  \begin{aligned} f(l,m) &= \xi_2 s^2 \\ g(l,m) &= \xi_2
    t^2 \end{aligned} \right\} = \{ \xi_2 s^2 = f
(q_1(\theta,\psi),q_2(\theta,\psi)) \}, \]
\[ C_3 = \left\{ 
  \begin{aligned} p_1(c,d) &= \xi_3  q_1(\theta,\psi) \\
    p_2(c,d) &= \xi_3 q_2(\theta,\psi) \end{aligned} \right\}, \qquad
C_4 = \left\{ \begin{aligned} F_2(x,y,u,v) &= 0 \\ G_2(x,y,u,v) &=
    0 \end{aligned} \right\}. \] Each of these curves is either a
double cover of $\PP^1$, or a quadric intersection in $\PP^3$.

The pushout functions $f_j \in K(C_{j-1})$ are as follows.  In each
case $C_j$ is $K$-birational to $\{ f_j (P) = \xi_j z^2 \} \subset
C_{j-1} \times \Gm$.  On $E = C_0$ the pushout function is $f_1 = x$.
On $C_1$ the pushout function is $f_2 = f(l,m)/m^2$.  We may solve by
linear algebra for $\alpha,\beta,\gamma \in K$ so that
\[ \alpha p_1(c,d) + \beta p_2(c,d) + \gamma p_3(c,d) = d^2. \]
Then $C_2$ has pushout function 
$f_3 = (\alpha q_1(\theta, \psi) + \beta q_2(\theta,\psi) 
 + \gamma s)/\psi^2$. Likewise we may solve for $r,s,t \in L$ so that
\[ r Q_1(\lambda,\mu) + s Q_2(\lambda,\mu) + t Q_3(\lambda,\mu) = \mu^2.\]
Let $\ell = r \theta + s \psi + t(c+ \gamma d)$ and
$\overline{\ell} = \overline{r} \theta + \overline{s} \psi + \overline{t}
(c+ {\overline{\gamma}} d)$, where the bars denotes the action of
$\Gal(L/K)$. Then $C_3$ has pushout function 
$f_4 = \ell \overline{\ell} /\psi^2$. 

The formulae in this section can be used to compute $S_m$ and $S'_m$
for $m \le 4$, by a method whose global part only requires that we
solve conics over $K$, and over quadratic extensions of $K$.  In
Section~\ref{sec:pushout} we extend to the case $m=5$, and in
Section~\ref{sec:res} we replace the conics over quadratic extensions
by quadric surfaces.

\section{Comparison with work of Bremner and Cassels}
\label{sec:compare}

The method used by Bremner and Cassels \cite{BC} to find rational
points of large height on $E': y^2 = x(x^2 + p)$, for $p$ a prime with
$p \equiv 5 \pmod{8}$, is a special case of the method outlined in
Section~\ref{sec:outline}. We now explain the relationship.  We have
tried to keep their notation, although in one place we switched $a,b$
to $s,t$, since we already used $a,b$ in (\ref{weqnE}).

We start with the elliptic curve $E : y^2 = x(x^2- 4p)$ over $K= \Q$.
The first stage of the argument in \cite{BC} implicitly switches to
this 2-isogenous curve. Then taking $\xi_1= p$ in~(\ref{cov1}) gives
$r^2 = p s^4 - 4 t^4$.  Since $p \equiv 5 \pmod{8}$ we can write $p=
u^2 + 4 v^2$ where $u,v$ are odd integers with $v \equiv 1
\pmod{4}$. A suitable parametrisation~(\ref{param1}) is now given by
\begin{align*}
f (l,m ) & = l^2 + m^2 \\
g (l,m ) & = v(l^2 - m^2) + u l m \\
h (l,m ) & = u(l^2 - m^2) - 4 v lm. 
\end{align*}
A local argument suggests taking $\xi_2 = 1$. Then~(\ref{cov2})
becomes $l^2 + m^2 = s^2$ and $v(l^2-m^2) + u lm = t^2$.  The first of
these conics is parametrised by $p_1 (c,d ) = c^2 - d^2$, $p_2 (c,d )
= 2cd$, $p_3 (c,d ) = c^2 + d^2$, and the second by
$q_1(\theta,\psi),q_2(\theta,\psi),q_3(\theta,\psi)$ (say).  If the
$q_i$ are suitably scaled then by local considerations it suffices to
take $\xi_3 = 1$. Then~(\ref{cov3}) becomes $q_1(\theta,\psi) = c^2 -
d^2$ and $q_2(\theta,\psi)=2cd$, equivalently
\begin{equation}
\label{conicQi}
 q_1(\theta,\psi) + i q_2(\theta,\psi) = (c + id)^2. 
\end{equation}
Parametrising a conic over $L= \Q(i)$ gives
\[ (\theta : \psi : c + i d ) = (Q_1(\lambda,\mu) : Q_2(\lambda,\mu) :
Q_3(\lambda,\mu) ). \] Then $\theta = \pi Q_1(\lambda,\mu)$ and $\psi
= \pi Q_2(\lambda,\mu)$ for some $\pi \in \Q(i)$.  A suitable value of
$\pi$ is determined by local considerations. Finally the procedure
described in Section~\ref{sec:outline} (with $\beta = i$) furnishes a
quadric intersection
\begin{equation}
\label{QIagain}
F_2(x,y,u,v) = G_2(x,y,u,v) = 0. 
\end{equation}
defined over $\Q$.

In summary, the examples of Bremner and Cassels are special in that,
(i) some of the conics can be parametrised ``for free'' and so do not
appear explicitly in their argument, (ii) the conics we do have to
solve are defined over either $\Q$ or $\Q(i)$, and (iii) by local
considerations only one choice of $\xi_j \in \Q^\times/(\Q^\times)^2$
need be considered at each stage ($j=1,2,3,4$).

Bremner and Cassels give a worked example in the case $p=877$.  We
record a few brief details. Firstly they take $u= 29$, $v= -3$ and
\begin{align*}
q_1(\theta,\psi) & =  -3 \theta^2 + 6 \theta \psi - 2 \psi^2 \\
q_2(\theta,\psi) & =  -7 \theta^2 - 4 \theta \psi - \psi^2 \\
q_3(\theta,\psi) & =  27 \theta^2 - 11 \theta \psi - 7 \psi^2.
\end{align*}
Then the conic~(\ref{conicQi}) is parametrised by
\begin{align*}
Q_1(\lambda,\mu) & = ( -4 - 3i)\lambda^2 
+ (10 - 22i)\lambda \mu + (29+6i) \mu^2  \\
Q_2(\lambda,\mu) & = ( -1 + 2i)\lambda^2 
+ (-16 - 6i)\lambda \mu + (6-29i) \mu^2  \\
Q_3(\lambda,\mu) & = ( -15 + 6i)\lambda^2 
+ (-46 - 70i)\lambda \mu + (76-75i) \mu^2.  
\end{align*}
By local considerations it suffices to take $\pi = 1+i$.  This leads
to a quadric intersection~(\ref{QIagain}).  Minimising and reducing,
as described in \cite{minred234}, suggests making the transformation
\begin{equation}
\label{minred}
\begin{aligned}
 x &=   x_1 + 2 x_2 - 4 x_3 - 6 x_4, &
 u &=    -x_1 - 2 x_3 + 2 x_4, \\
 y &=  2 x_1 - x_2 + 6 x_3 - 4 x_4, &
 v &=    x_2 - 2 x_3 - 2 x_4
\end{aligned}
\end{equation}
whereupon~(\ref{QIagain}) becomes
\begin{equation}
\label{4cov}
\begin{aligned}
 x_1 x_2 + x_1 x_3 + x_1 x_4 - x_2 x_3 + x_2 x_4 
       + x_3^2 - 2 x_3 x_4 - x_4^2  &= 0 \\
 x_1^2 - 4 x_1 x_2 + 3 x_1 x_3 + x_1 x_4 - x_2^2 
       - x_2 x_3 + 3 x_2 x_4 + 4 x_3 x_4 &= 0. 
\end{aligned}
\end{equation}
A little searching finds the rational points
\[(x_1 : x_2 : x_3 : x_4 ) = (-2 : 57 : 85 : 16), \, (57 : 2 : -16 :
85). \] We substitute in~(\ref{minred}) to recover the solutions
\[(x : y : u : v ) = (324 : -385 : 136 : 145), \, (385 : 324 : -145 :
136) \] found by Bremner and Cassels. These points map down to the
same point of infinite order $P=(x_P,y_P)$ on $E(\Q)$.  The
co-ordinates of $P$ are
\begin{align*}
x_P &= -292214148680270491236/4612160965^2 \\
y_P &= 20949922565086352416107761007588/4612160965^3. 
\end{align*}
The point recorded in \cite{BC} is the image of $P$ under the
2-isogeny $\phi : E \to E'$, and accordingly has (canonical) height
twice that of $P$.

The general implementation of 4-descent in Magma (due to
Womack~\cite{Womack} and Watkins) is able to find the
4-covering~(\ref{4cov}) in a couple of seconds. However the method it
uses involves computing the class group and units of a degree 4 number
field; in this case $\Q(\sqrt{6 - 29i})$.  In contrast the method of
Bremner and Cassels only requires that we solve conics over $\Q$ and
$\Q(i)$.

\section{Computing pushout forms}
\label{sec:pushout}

Let $C \subset \PP^3$ be a non-singular quadric intersection.  The $4$
singular fibres in the pencil of quadrics defining $C$ are cones over
conics $\Gamma_1, \ldots ,\Gamma_4$. Projecting from the vertex of
each cone gives a degree-$2$ morphism $\nu_i : C \to \Gamma_i$ with
fibre $\aaa_i$. By considering the tangent plane at a smooth point on
the cone over $\Gamma_i$ we see that $2 \aaa_i \sim H$ where $H$ is
the hyperplane section on $C$. Therefore the differences $\aaa_i -
\aaa_j$ represent elements of order $2$ in $\Pic^0(C)$.

Let $E$ be the Jacobian of $C$. By the previous paragraph $\{
\Gamma_1, \ldots, \Gamma_4 \}$ is a torsor under $E[2]$.  In
particular, there is a Galois equivariant bijection between $E[2]
\setminus \{0\}$ and the partitions of the singular fibres into $2$
sets of $2$. We are interested in elliptic curves with a rational
$2$-torsion point. We therefore fix $0 \not= T \in E(K)[2]$, and order
the $\Gamma_i$ so that $T$ corresponds to the partition $(\Gamma_1,
\Gamma_2; \Gamma_3, \Gamma_4)$. In other words $\aaa_1 - \aaa_2 \sim
\aaa_3 - \aaa_4$ represents the class of $T$ in $\Pic^0(C) \isom E$.

Let $\rho : \Gal(\Kbar/K) \to S_4$ describe the action of Galois on
the $\Gamma_i$. Since $T$ is $K$-rational we have $\im(\rho) \subset
\langle (1324),(12) \rangle$. We now make the assumption that
\begin{equation}
\label{galois-hypothesis}
 \im(\rho) \subset \langle (12),(34) \rangle. 
\end{equation}
We distinguish two possibilities. 
\begin{itemize}
\item (split case) The conics $\Gamma_1$ and $\Gamma_2$ are defined
over $K$.
\item (non-split case) The conics $\Gamma_1$ and $\Gamma_2$ are
  defined over a quadratic extension $L/K$, and are
  $\Gal(L/K)$-conjugates.
\end{itemize}
Let $\Gamma = \Gamma_1 \times \Gamma_2$, respectively $\Res_{L/K}
\Gamma_1$.  If $C/K$ is everywhere locally soluble then by the Hasse
Principle for conics, we have $\Gamma(K) \not= \emptyset$. A point $P
\in \Gamma(K)$ corresponds to a pair of points $P_1 \in \Gamma_1$ and
$P_2 \in \Gamma_2$ which are either defined over $K$ or are
$\Gal(L/K)$-conjugates. The tangent planes to $P_1$ and $P_2$, as
points on the cones over $\Gamma_1$ and $\Gamma_2$, are defined by
linear forms $\ell_1$ and $\ell_2$. In the split case these are
defined over $K$.  In the non-split case we may arrange that they are
$\Gal(L/K)$-conjugates.

We say that a quadratic form $f \in K[x_1,\ldots,x_4]$ is a 
{\em pushout form} on $C$ if $f/x_1^2 \in K(C)$ is a pushout function,
i.e.  $\divv(f/x_1^2) = 2 \bbb$ for some divisor $\bbb$ representing
the class of $T$ in $\Pic^0(C) \isom E$.

\begin{Lemma}
\label{lem1}
$f = \ell_1 \ell_2 \in K[x_1,\ldots,x_4]$ is a pushout form on $C$.
\end{Lemma}
\begin{Proof}
We have $\divv(f/x_1^2) = 2(\aaa_1 + \aaa_2 - H)$ and 
$\aaa_1 + \aaa_2 - H \sim \aaa_1 - \aaa_2$. 
\end{Proof}

Let $f = \ell_1 \ell_2$ as above and let $\pi : D \to C$ be the
degree-$2$ covering with $K(D) = K(C)(\sqrt{f})$. In other words $D$
is the smooth curve of genus one $K$-birational to $\{ f(P) = z^2 \}
\subset C \times \Gm$. We now show how to write $D$ as a quadric
intersection with hyperplane section $\pi^* \aaa_1 \sim \pi^*\aaa_2$.

We suppose we are in the non-split case, the split case being similar.
The conic $\Gamma_1$ is defined by a quadratic form $q \in L[X,Y,Z]$,
and the cone above $\Gamma_1$ has equation $q(m_1,m_2,m_3) = 0$ for
some linear forms $m_1,m_2,m_3 \in L[x_1,\ldots,x_4]$. Since we have
found an $L$-rational point on $\Gamma_1$, we may parametrise this
conic, say
\[(X:Y:Z) = (Q_1(\lambda,\mu) : Q_2(\lambda,\mu) : Q_3(\lambda,\mu))\]
where $Q_1,Q_2$ and $Q_3$ are binary quadratic forms defined over $L$.
If $L = K(\beta)$ this gives equations 
\begin{equation}
\label{three-eqns}
m_i (x_1, \ldots, x_4) = Q_i(x+\beta y,u+\beta v) 
\end{equation}
for $i=1,2,3$. Writing each side in terms of the basis $1,\beta$ for
$L$ over $K$ we get $6$ equations with coefficients in $K$.  The left
hand sides are linear forms in $x_1, \ldots, x_4$ and the right hand
sides are quadratic forms in $x,y,u,v$. By taking linear combinations
we obtain equations
\[   F(x,y,u,v) = G(x,y,u,v) = 0 \]
and $x_i = r_i(x,y,u,v)$ for $i=1,\ldots,4$. 

\begin{Theorem}
\label{thm2}
The curve $D = \{F = G = 0\} \subset \PP^3$ is a non-singular quadric
intersection, and the map $\pi = (r_1: \ldots:r_4):D \to C$ is a
morphism of degree $2$. Moreover $D$ has hyperplane section $\pi^*
\aaa_1 \sim \pi^* \aaa_2$.
\end{Theorem}
\begin{Proof}
  For the proof we are free to extend our field $K$ and to make
  changes of co-ordinates. So we may suppose
\begin{equation}
\label{geomC}
C = \left\{ 
\begin{aligned} f(l,m) &= s^2 \\ g(l,m) &= t^2 \end{aligned}
\right\} \subset \PP^3. 
\end{equation}
We parametrise the conics $f(l,m) = s^2$ and $g(l,m) = t^2$ as
\begin{equation}
\label{conicparam}
 \begin{aligned}
( l : m : s ) & = ( p_1(c,d) : p_2(c,d) : p_3(c,d) ), \\
( l : m : t ) & = ( q_1(\theta,\psi) : q_2(\theta,\psi) : q_3(\theta,\psi) ). 
\end{aligned}
\end{equation}
Then by the above construction (in the split case) we have
\begin{equation}
\label{geomD}
D = \left\{ 
\begin{aligned} p_1(c,d) &= q_1(\theta,\psi) \\ 
p_2(c,d) &= q_2(\theta,\psi) \end{aligned} \right\} \subset \PP^3. 
\end{equation}

If a quadric intersection is given by $4 \times 4$ symmetric matrices
$A$ and $B$, then we may associate to it the binary quartic $F(x,y) =
\det (A x + B y)$. To prove that the quadric intersection is
non-singular it suffices to show that $F$ has distinct roots in
$\PP^1$. We write $\Delta$ for the discriminant of a binary quadratic
form. Since the linear combinations of $p_1$ and $p_2$ that are
perfect squares can be computed from the roots of $f$, we have
\begin{align*} 
  \Delta( m p_1 - l p_2) = \kappa f(l,m) \\
  \Delta( m q_1 - l q_2) = \kappa' g(l,m) 
\end{align*}
for some $\kappa,\kappa' \in K^\times$. 
Therefore the binary quartic associated to 
$D$ is a scalar multiple of $f(l,m)g(l,m)$. Since $C$ defined
by~(\ref{geomC}) is non-singular, it now follows that $D$ 
defined by~(\ref{geomD}) is non-singular.

The morphism $\pi : D \to C$ is given by
\[ (l:m:s:t) = (p_1(c,d):p_2(c,d):p_3(c,d):q_3(\theta,\psi)). \] It is
easy to see that $\pi$ has degree $2$ with fibres of the form
$\{(c:d:\theta:\psi),(c:d:-\theta:-\psi)\}$. Moreover if $\aaa_1 =
(l:m:s:t) + (l:m:s:-t)$ then $\pi^* \aaa_1$ is the hyperplane section
given by solving the first equation in~(\ref{conicparam}) for $(c:d)$.
Likewise if $\aaa_2 = (l:m:s:t) + (l:m:-s:t)$ then $\pi^* \aaa_2$ is
the hyperplane section given by solving the second equation
in~(\ref{conicparam}) for $(\theta:\psi)$.
\end{Proof}

Let $\phi : E \to E'$ be the $2$-isogeny with kernel $\{ 0,T \}$. 
Let $\phihat : E' \to E$ be the dual isogeny, say with kernel
$\{0,T'\}$. 

\begin{Theorem} 
\label{thm3}
The degree-$2$ covering $\pi : D \to C$ constructed in 
Theorem~\ref{thm2} is a $\phihat$-covering, i.e. there is a 
commutative diagram
\begin{equation*}
\xymatrix{ D \ar[d]_{\isom} \ar[r]^{\pi} & C \ar[d]_{\isom} \\
E' \ar[r]^\phihat & E }
\end{equation*}
where the vertical maps are isomorphisms defined over $\Kbar$.
\end{Theorem}
\begin{Proof}
We give details in the non-split case. Recall that
$\Gamma_1$ is parametrised
by binary quadratic forms $Q_1,Q_2,Q_3$. We solve for $r,s,t$
so that 
$ r Q_1(\lambda,\mu) + s Q_2(\lambda,\mu) + t Q_3(\lambda,\mu) = \mu^2$.
We then take $f = \ell_1 \ell_2$ where
$\ell_1 = r m_1 + s m_2 + t m_3$ and $\ell_2$ is its $\Gal(L/K)$-conjugate.
By~(\ref{three-eqns}) we have
\begin{align*}
\ell_1(x_1,\ldots,x_4) &= (u + \beta v)^2 \\
\ell_2(x_1,\ldots,x_4) &= (u + \overline{\beta} v)^2 
\end{align*}
Then $t = (u + \beta v)(u + \overline{\beta} v) \in K[x,y,u,v]$ is 
a quadratic form with
\begin{equation}
\label{eqn:odd}
 f(r_1, \ldots, r_4) \equiv t^2 \mod{I(D)}. 
\end{equation}
Therefore $\pi^*(f/x_1^2)$ is the square of a rational function in $K(D)$.
Since $f$ is a pushout form (corresponding to $T$), 
the result follows.
\end{Proof}

The pushout form $f \in K[x_1,\ldots,x_4]$ and quadric intersection $D
\subset \PP^3$ were constructed from a pair of points on $\Gamma_1$
and $\Gamma_2$. We may equally construct a pushout form $f^\dagger \in
K[x_1,\ldots,x_4]$ and quadric intersection $D^\dagger \subset \PP^3$
from a pair of points on $\Gamma_3$ and $\Gamma_4$. We have
$\divv(f/x_1^2) = 2 \bbb$ and $\divv(f^\dagger/x_1^2) = 2
\bbb^\dagger$ with $\bbb \sim \bbb^\dagger$ representing the class of
$T$ in $\Pic^0(C) \isom E$. Therefore $f/f^\dagger = c h^2$ for some
$c \in K$ and $h \in K(C)$. Scaling $f$ and $f^\dagger$ appropriately
we may assume $c=1$. The quadric intersections $D$ and $D^\dagger$ are
now isomorphic as curves (indeed as $\phihat$-coverings of $C$) but
have different hyperplane sections $H$ and $H^\dagger$.

\enlargethispage{0.8ex}

\begin{Theorem}
\label{thm4}
The divisor $H - H^\dagger$ represents the class of $T'$ in 
$\Pic^0(D) \isom E'$. 
\end{Theorem}
\begin{Proof}
Let $\pi : D \to C$ be the degree-$2$ covering as above.
By Theorem~\ref{thm2} we have 
\[ H \sim \pi^* \aaa_1 \sim \pi^* \aaa_2 \quad \text{ and } \quad
H^\dagger \sim \pi^* \aaa_3 \sim \pi^* \aaa_4. \] Let $E[2] =
\{0,T,S_1,S_2\}$ with $\aaa_1 - \aaa_3$ representing the class of
$S_1$ in $\Pic^0(C) \isom E$. By Theorem~\ref{thm3} there is a
commutative diagram
\begin{equation*}
  \xymatrix{ \Pic^0(C) \ar[d]_{\isom} \ar[r]^{\pi^*} & \Pic^0(D) \ar[d]_{\isom} \\
    E \ar[r]^\phi & E' }
\end{equation*}
Then $H - H^\dagger \sim \pi^*(\aaa_1 - \aaa_3)$ represents the class
of $\phi S_1 = T'$ in $\Pic^0(D) \isom E'$.
\end{Proof}

We will need explicit equations for the isomorphism between the
quadric intersections $D$ and $D^\dagger$.  We compute these as
follows. The $2$-uple embedding $D \subset \PP^7$ is defined by $8$
quadratic forms, chosen so that together with the equations for $D$
they give a basis for the space of all quadratic forms on
$\PP^3$. Since $2 H \sim 2 H^\dagger$ the $2$-uple embeddings of $D$
and $D^\dagger$ are related by a change of co-ordinates on $\PP^7$. We
now explain how to choose $8$ quadratic forms for $D$ in a particular
way, so that when we repeat for $D^\dagger$, the change of
co-ordinates needed on $\PP^7$ is trivial.

Since $\pi : D \to C$ is a degree-$2$ Galois covering, the $8$
quadratic forms may be given as $4$ ``even'' forms $r_1,\ldots, r_4$
and $4$ ``odd'' forms $s_1, \ldots, s_4$.  The even forms $r_1,
\ldots, r_4$ are those giving the morphism $\pi : D \to C$ in
Theorem~\ref{thm2}. The quadratic form $t$ computed in~(\ref{eqn:odd})
is an example of an odd form. Further odd forms may be computed as
follows.

Once we have found one pushout form $f$ on $C$, it is easy to find
more by computing Riemann-Roch spaces. (With the refinements in
Section~\ref{sec:res} it turns out that we have this information
anyway.)  Let $f_i$ be another pushout form, scaled so that $f f_i
\equiv h_i^2 \mod{I(C)}$ for some quadratic form $h_i$.  Since
$\pi^*(f_i /x_1^2)$ is a square in $K(D)$, and $D$ is projectively
normal, there is a quadratic form $s_i$ such that
\begin{equation}
\label{eqn:odd2}
 f_i(r_1, \ldots,r_4) \equiv s_i^2 \mod{I(D)}.
\end{equation}
It follows by~(\ref{eqn:odd}) and~(\ref{eqn:odd2}) that if we scale
the forms $s_i$ appropriately then
\[ h_i(r_1,\ldots,r_4) \equiv s_i t \mod{I(D)}. \] We use this last
equation to solve for $s_i$ by linear algebra. Repeating for pushout
forms $f_1, \ldots, f_4$ gives odd forms $s_1, \ldots, s_4$ are
required.

We now have quadric intersections $D \subset \PP^3$ and $D^\dagger
\subset \PP^3$ and an isomorphism between their $2$-uple embeddings,
given by a change of co-ordinates on $\PP^7$. We say that $D$ and
$D^\dagger$ are {\em companion quadric
  intersections}. Theorem~\ref{thm4} shows that, under this change of
co-ordinates, the square of a linear form on $D^\dagger$ corresponds
to a pushout form on $D$.  This is exactly what we need for the next
stage of the descent.

The results of this section may be applied as follows.  Let $C_3$ be
the $2 \phihat$-covering of $E$ considered in
Section~\ref{sec:outline}. Then $C_3$ satisfies the Galois
hypothesis~(\ref{galois-hypothesis}). So we may construct a
$\phi$-covering of $C_3$ (and hence a $4$-covering over $E$) in the
form of companion quadric intersections $C_4$ and $C^\dagger_4$. The
construction of $C_4$ was already described in
Section~\ref{sec:outline}. However the advantage of also computing
$C^\dagger_4$, is that the square of a linear form on $C^\dagger_4$
corresponds to a pushout form on $C_4$.  We can then compute the
pairing $\Theta_4 : S_4 \times S'_4 \to \F_2$ in
Theorem~\ref{thm:pair}. In conclusion, we can compute $S_m$ and $S'_m$
for $m \le 5$ by a method whose global part only involves solving
conics over $K$, and over quadratic extensions of $K$.

\section{Extension to $8$-descent}
\label{sec:8}

Suppose that $E$ has full rational $2$-torsion, say 
\[E(K)[2] = \{0,T_1,T_2,T_3\}.\] Let $C \subset \PP^3$
be an everywhere locally soluble $4$-covering of $E$. 
Repeating the method of Section~\ref{sec:pushout} three times,
gives pushout forms $f_i \in K[x_1,\ldots,x_4]$ on $C$ with
$\divv(f_i/x_1^2) = 2 \bbb_i$ where $\bbb_i$ represents the class
of $T_i$ in $\Pic^0(C) \isom E$. The rational functions $f_i/x_1^2$
are now exactly those required to compute the Cassels-Tate
pairing
\[ S^{(4)}(E/K) \times S^{(2)}(E/K) \to \Q/\Z \] following
Swinnerton-Dyer's generalisation~\cite{SD} of the method of
Cassels~\cite{Ca98}.  This allows us to compute the pairing $\Theta'_5
: S'_5 \times S'_5 \to \F_2$ and hence its kernel $S'_6$.  By
Lemma~\ref{lem:pn} the upper bound
\[ \rank E(K) \le \dim S_5 + \dim S'_6 -2 \]
is the same as that obtained by $8$-descent on $E$.

The following lemma, which we state and prove using the notation and
conventions of Section~\ref{sec:book}, describes the necessary
bookkeeping.

\begin{Lemma} 
\label{lem:book8}
Let $\alpha_1, \beta_1 \in S_5$. Suppose $\alpha_4,
\beta_4 \in \Sel(4)$ and $\beta_2 \in \Sel(2)$ with $\alpha_4 \mapsto
\alpha_1$ and $\beta_4 \mapsto \beta_2 \mapsto \beta_1$. Then
\begin{enumerate}
\item There exists $\xi \in S_3$ with 
$\Theta_3(\xi,\eta) + \langle \beta_4, \eta \rangle = 0$
for all $\eta \in S_3$.
\item For any $\xi$ satisfying (i) we have
\begin{equation}
\label{eq:th5}
 \Theta_5(\alpha_1, \beta_1) = \langle \alpha_4, \xi + \beta_2 \rangle. 
\end{equation}
\end{enumerate}
\end{Lemma}

\begin{Proof}
  Let $\beta_5 \in \Sel(5)$ with $\beta_5 \mapsto \beta_1$. The images
  of $\beta_5$ in $\Sel(4)$ and $\Sel(2)$ are $\xi' + \beta_4$ and
  $\xi + \beta_2$ for some $\xi' \in \Sel(3)$ and $\xi \in \Sel(1)$
  with $\xi' \mapsto \xi$. Then $\xi \in S_3$ and for any $\eta \in
  S_3$ we have
  \[ \Theta_3(\xi,\eta) + \langle \beta_4, \eta \rangle = \langle \xi'
  + \beta_4, \eta \rangle = 0. \] Moreover $\Theta_5(\alpha_1,\beta_1)
  = \langle \alpha_1, \beta_5 \rangle = \langle \alpha_4, \xi +
  \beta_2 \rangle$. To complete the proof we must show that if we
  replace $\xi$ by $\xi + \gamma_1$ for some $\gamma_1 \in S_4$ then
  the formula~(\ref{eq:th5}) is unchanged.  However if $\gamma_4
  \in \Sel(4)$ with $\gamma_4 \mapsto \gamma_1$ then $\langle
  \alpha_4, \gamma_1 \rangle = \langle \alpha_1, \gamma_4 \rangle =0$
  where for the second equality we use that $\alpha_1 \in S_5$.
\end{Proof}

\section{Restriction of scalars}
\label{sec:res}

Our methods so far rely on being able to solve conics over $K$, and
over quadratic extensions of $K$. In this section we show how to
replace the latter problem with that of solving a quadric surface over
$K$. This is to our advantage since, in the case $K=\Q$, there is a
particularly efficient method due to D. Simon \cite{SimonQuadrics} for
solving rank $4$ quadratic forms.

Let $\Gamma_1$ be a conic defined over a quadratic extension
$L/K$. The restriction of scalars $\Gamma = \Res_{L/K} \Gamma_1$ is a
degree $8$ del Pezzo surface in $\Res_{L/K} \PP^2 \subset \PP^8$.
Alternatively, by passing to an affine piece first, one gets that
$\Gamma$ is birational to an intersection of two quadrics in $\PP^4$
with $4$ singular points at infinity.  We show however that for the
conics arising in our descent calculations, we can write $\Gamma$ as a
quadric surface in $\PP^3$.

Let $E/K$ be an elliptic curve with a fixed rational $2$-torsion point
$T$. Suppose $C$ and $C^\dagger$ are companion quadric intersections
with respect to $T$. By this we mean that $C$ and $C^\dagger$ are
isomorphic as curves, but the difference of their hyperplane sections
represents $T$ in $\Pic^0(C) \isom E$. Let $\Gamma_1, \ldots,
\Gamma_4$ be the conics associated to $C$, ordered as specified at the
start of Section~\ref{sec:pushout}. Suppose we are in the non-split
case, i.e. $\Gamma_1$ and $\Gamma_2$ are $\Gal(L/K)$-conjugates. Let
$\Gamma = \Res_{L/K} \Gamma_1$.  As described in Lemma~\ref{lem1},
each $P \in \Gamma$ gives rise to a pushout form on $C$. Since $C$ and
$C^\dagger$ are companion quadric intersections, this in turn
corresponds to a linear form on $C^\dagger$. There is therefore a
natural map $\Gamma \to (\PP^3)^\vee$ where $(\PP^3)^\vee$ is the dual
of the ambient space for $C^\dagger$.

If $\{ Q = 0 \} \subset \PP^3$ is a non-singular quadric surface then
mapping each point to its tangent plane gives the dual quadric surface
$\{ Q' = 0 \} \subset (\PP^3)^\vee$. The symmetric matrix representing
$Q'$ is the inverse of that representing $Q$.

\begin{Theorem}
\label{thm5}
The image of $\Gamma \to (\PP^3)^\vee$ is a non-singular quadric
surface, and is dual to one of the quadrics in the pencil defining
$C^\dagger$.
\end{Theorem}
\begin{Proof}
  For the proof we are free to extend our field $K$ and make changes
  of co-ordinates. We may therefore suppose that $C$ and $C^\dagger$
  are the images of $E : y^2 = x(x^2+ ax+b)$ under the linear systems
  $|4.0|$ and $|3.0 + T|$ where $T = (0,0)$. Explicitly
\begin{equation*}
C = \left\{ 
\begin{aligned} x_1 x_4 &= x_2^2 \\ x_3^2 &= 
 x_2(x_4 + a x_2 + b x_1) \end{aligned}
\right\} \subset \PP^3
\end{equation*}
is the image of $E$ under $(x_1:x_2:x_3:x_4) = (1:x:y:x^2)$, and
\begin{equation}
\label{qi-cc}
C^\dagger = \left\{ 
\begin{aligned} x_1 x_3 &= x_2 x_4 \\ x_3 x_4 &= 
 x_2^2 + a x_1 x_2 + b x_1^2 \end{aligned}
\right\} \subset \PP^3 
\end{equation}
is the image of $E$ under $(x_1:x_2:x_3:x_4) = (1:x:y:y/x)$.

Let $P_1 = (\alpha_1 : \alpha_2: 0 : \alpha_4)$ and
$P_2 = (0: \beta_2 : \beta_3 : \beta_4)$ be points on the rank~$3$ 
quadrics defining $C$. The tangent planes at these points
are defined by linear forms 
\begin{align*}
\ell_1 &= \alpha_4 x_1 - 2 \alpha_2 x_2 + \alpha_1 x_4, \\
\ell_2 &= (\beta_4 + a \beta_2) x_2 + \beta_2 (x_4 + a x_2 + b x_1)
 - 2 \beta_3 x_3.
\end{align*}
Using the relations $\alpha_1 \alpha_4 = \alpha_2^2$ and
$\beta_3^2 = \beta_2(\beta_4 + a \beta_2)$ we find that
\[ \alpha_4 \beta_2 \ell_1(1,x,y,x^2) \ell_2(1,x,y,x^2) =
m(1,x,y,y/x)^2 x \] in $K(E)$, where
\[ m(x_1,\ldots,x_4) = \alpha_4 \beta_3 x_1 - \alpha_2 \beta_3 x_2
 + \alpha_2 \beta_2 x_3 - \alpha_4 \beta_2 x_4. \]
The map $\Gamma \to (\PP^3)^\vee$ is therefore given by
\[ (P_1,P_2) \mapsto ( \alpha_4 \beta_3 : - \alpha_2 \beta_3 :
 \alpha_2 \beta_2 : - \alpha_4 \beta_2). \]
The image of this map has equation $y_1 y_3 = y_2 y_4$. This is dual to
the first of the quadrics in~(\ref{qi-cc}). 
\end{Proof} 

\begin{Remark}
  The binary quartic associated to $C^\dagger$ defines a double cover
  of $\PP^1$, again with Jacobian $E$. Translation by the $2$-torsion
  point $T$ induces an involution of $\PP^1$. The quadric in the
  pencil defining $C^\dagger$ arising in Theorem~\ref{thm5}
  corresponds to one of the fixed points of this involution. (The
  other fixed point corresponds to the quadric that arises when we try
  to solve $\Gamma_3$ and $\Gamma_4$.)  This may be seen by inspection
  of the above proof. Indeed the binary quartic associated
  to~(\ref{qi-cc}) is
  \[F(u,v) = u^4 + 2 a u^2 v^2 + (a^2 - 4b) v^4 \] and the involution
  induced by translation by $T$ is $(u:v) \mapsto (u:-v)$. The fixed
  point $(u:v) = (1:0)$ then corresponds to the first of the quadrics
  in~(\ref{qi-cc}).
\end{Remark}

Solving the quadric surface in Theorem~\ref{thm5} gives a linear form
on $C^\dagger$. Suppose we know the change of co-ordinates relating
the $2$-uple embeddings of $C$ and $C^\dagger$. This then converts the
square of our linear form on $C^\dagger$ to a pushout form $f$
on~$C$. We write $C \cap \{f = 0\} = 2 \bbb$ where $\bbb$ is a degree
$4$ effective divisor. By construction, the push-forward of $\bbb$ via
$\nu_1 : C \to \Gamma_1$ contains an $L$-rational point in its
support.  We can then solve for $P_1 \in \Gamma_1(L)$ as required.

To make use of the above refinements, we must arrange that at each
stage of the descent, our covering curve is represented by a pair of
companion quadric intersections, and that we know the change of
co-ordinates on $\PP^7$ relating their $2$-uple embeddings.

Let $C : y^2 = f(l,m)g(l,m)$ be the $\phihat$-covering of $E$ we
called $C_1$ in Section~\ref{sec:outline}. Then $C$ is a double cover
of $\PP^1$ with fibre $F$ and ramification points $P_1, \ldots, P_4$,
say with $P_1,P_2$ corresponding to the roots of $f$, and $P_3,P_4$
corresponding to the roots of $g$. We have $2P_i \sim F$ and $P_1 +
P_2 + P_3 + P_4 \sim 2 F$. Moreover $P_1- P_2 \sim P_3 - P_4$
represents the class of $T'$ in $\Pic^0(C) \isom E'$.  We now let
$C_1$ and $C_1^\dagger$ be the images of $C$ under the linear systems
$|2F|$ and $|F + P_1 + P_2|$. Explicitly
\begin{equation*}
C_1 = \left\{ 
\begin{aligned} x_1 x_3 &= x_2^2 \\ x_4^2 &= 
 \xi_1 x_1^2 + a x_1 x_3 + (b/\xi_1) x_3^2 \end{aligned}
\right\} \subset \PP^3 
\end{equation*}
is the image of $C$ under $(x_1:x_2:x_3:x_4) = (f(l,m):y:g(l,m):h(l,m))$, and
\begin{equation*}
C_1^\dagger = \left\{ 
\begin{aligned} x_1 x_4 &= x_2 x_3 \\ f(x_1,x_2) &= 
 g(x_3,x_4) \end{aligned}
\right\} \subset \PP^3 
\end{equation*}
is the image of $C$ under 
\begin{align*}
(x_1:x_2:x_3:x_4) &= (yl:ym:f(l,m)l:f(l,m)m) \\
&= (g(l,m)l:g(l,m)m:yl:ym). 
\end{align*}
It is straightforward to work out the change of co-ordinates relating
the $2$-uple embeddings of $C_1$ and $C_1^\dagger$.

Starting with $C_1$ and $C_1^\dagger$ we apply our method to compute a
$\phi$-covering of $C_1$ (and hence a $2$-covering of $E$) in the form
of companion quadric intersections $C_2$ and $C_2^\dagger$. The curve
$C_2$ is the same as that in Section~\ref{sec:outline}, and as there
it is computed by solving conics over $K$. To compute $C_2^\dagger$ we
must solve the quadric surface $f(x_1,x_2) = g(x_3,x_4)$ that appears
as the second equation for $C_1^\dagger$. Since we already solved the
conics $f(l,m) = \xi_2 s^2$ and $g(l,m) = \xi_2 t^2$ in
Section~\ref{sec:outline}, we can read off a solution ``for free''.

At the next stage we compute a $\phihat$-covering of $C_2$ (and hence
a $2 \phihat$-covering of $E$) in the form of companion quadric
intersections $C_3$ and $C_3^\dagger$.  The curve $C_3$ is the same as
that in Section~\ref{sec:outline}, and as there it is computed by
solving conics over $K$. However to compute $C_3^\dagger$ we must
solve a quadric surface over $K$. In the final stage, we compute a
$\phi$-covering of $C_3$ (and hence a $4$-covering of $E$) in the form
of companion quadric intersections $C_4$ and $C_4^\dagger$.  To
compute each of these we must solve a quadric surface over $K$.

This is the limit of our method since $C_4$ does not satisfy the
Galois hypothesis~(\ref{galois-hypothesis}). Indeed the singular
fibres in the pencil defining $C_4$ are in general defined over a
degree $4$ extension of $K$. Solving a conic over a degree $4$ number
field is not practical for the examples we have in mind.  However if
we could solve such conics, then we could compute a $2$-covering of
$C_4$ (and hence an $8$-covering on $E$) as described in
\cite{Stamminger}.

\section{Examples}
\label{sec:examples}

We have written a program in~Magma \cite{magma} for performing the
higher descents described in this paper, for elliptic curves over $\Q$
with a rational $2$-torsion point.  The Magma functions we use for
solving quadratic forms of ranks 3 and 4 over $\Q$ are based on
\cite{SimonConics}, \cite{SimonQuadrics}.  See also \cite{CR} for the
rank 3 case.  Our program is available in Magma (version 2.21) as the
function {\tt TwoPowerIsogenyDescentRankBound}.

We have used our program to help search for curves of large rank in
certain families of elliptic curves.  For example, we found the first
example of an elliptic curve $E/\Q$ with $E(\Q) \isom \Z/12\Z \times
\Z^4$.  We also found $10$ new examples of elliptic curves $E/\Q$ with
$E(\Q) \isom \Z/2\Z \times \Z/8\Z \times \Z^3$, including one where
every point of infinite order has canonical height greater than
$100$. These examples are listed on Dujella's website \cite{Dujella}.
We have also contributed to a project run by Mark Watkins
\cite{Watkins+}, searching for congruent number elliptic curves of
large rank.

The main reason we need higher descents for these searches is that we
would otherwise be swamped by examples which, while appearing to be
candidates for large rank, turn out instead to have large $2$-primary
part of $\Sha$.  The examples we have chosen to present in this
section are therefore elliptic curves which our program was quickly
able to show do {\em not} have large rank.

We do not give details of every step in computing the covering curves
and pushout functions. However the answers may be checked as
follows. Each of our covering curves is either a double cover of
$\PP^1$, with equation $y^2 = g(x,z)$ where $g$ is a binary quartic,
or a quadric intersection in $\PP^3$. In both cases classical
invariant theory gives a formula for the Jacobian.  In the second
case, we may represent the quadric intersection $C_4 \subset \PP^3$ by
a pair of $4 \times 4$ symmetric matrices $A$ and $B$. Then $C_4$ is a
2-covering of $y^2 = g(x,z)$ where $g(x,y) = \det (A x + B y)$. At
each stage our program makes changes of co-ordinates to simplify the
equations for these covering curves, using the algorithms in
\cite{minred234}.  In checking the examples, it is useful to note that
there are algorithms implemented in Magma for testing equivalence of
binary quartics (see \cite{testeqbq}) and quadric intersections (see
\cite{improve4}).

The pushout functions on $y^2 = g(x,z)$ take the form $f = (y -
\lambda_2(x,z))/z^2$ where $\lambda_2$ is a binary quadratic form. The
pushout functions on a quadric intersection take the form $f =
\lambda_4(x_1, \ldots, x_4)/x_1^2$ where $\lambda_4$ is a quadratic
form. We say that $y - \lambda_2(x,z)$ and $\lambda_4(x_1,\ldots,
x_4)$ are {\em pushout forms}. There are several ways of checking that
$\divv(f) = 2 \bbb$ for some divisor $\bbb$. For example in Magma one
could compute resultants, or use Groebner bases, or use the function
field machinery. If $f$ is not a constant times the square of a
rational function, and the Jacobian has only one rational $2$-torsion
point $T$, then such a calculation proves that $f$ is a pushout
function corresponding to $T$.

\begin{Example}
Let $E/\Q$ be the elliptic curve $y^2 = x(x^2 + a x + b)$ where
\begin{align*}
a &= 91502230365284038, \\
b &= 489792722057841784540058275212361.
\end{align*}
This is an example with $E(\Q)_{\tors} \isom \Z/12\Z$.
There is a $2$-isogeny $\phi : E \to E'$ where 
$E'$ has equation $y^2 = x(x^2 + a' x + b')$
and $a' = -2a$, $b' = a^2 - 4b$.  With notation as in Section~\ref{intro},
we find that 
\begin{align*}
S_1 = S_2 = S_3 &= \langle 15, 73, 87, 231, 28619 \rangle \subset \Q^\times/(\Q^\times)^2, \\
S'_1 = S'_2 = S'_3 &= \langle -272196179 \rangle \subset \Q^\times/(\Q^\times)^2. 
\end{align*}
Therefore $\rank E(\Q) \le 4$. To improve this upper bound for 
the rank we compute the pairing $\Theta_3 : S_3 \times S_3 \to \F_2$.

Let $\xi_1 = 15 \times 44660^2$. Parametrising a conic, we find binary 
quadratic forms
\begin{align*}
f(l,m) &=    195346817865 l^2 - 490516840068 l m + 33576198052 m^2, \\
g(l,m) &=   -1473071 l^2 + 7386682 l m - 9255011 m^2, \\
h(l,m) &=    348523798338106 (2067 l^2 - 10091 l m + 12314 m^2)
\end{align*}
satisfying $ \xi_1 f^2 + a' f g + (b'/\xi_1) g^2 = h^2$.
Parametrising two further conics, gives binary quadratic forms
\begin{align*}
  p_1(c,d) &=  65928 c^2 + 582550 c d + 1159554 d^2, \\
  p_2(c,d) &=  13590 c^2 + 429375 c d - 202059 d^2, \\
  p_3(c,d) &=  8585676 (2375 c^2 - 6774 c d - 71700 d^2)
\end{align*}
and
\begin{align*}
  q_1(\theta,\psi) &=  819 \theta^2 + 1717 \theta \psi + 3725 \psi^2, \\
  q_2(\theta,\psi) &=  329 \theta^2 + 717 \theta \psi + 1510 \psi^2, \\
  q_3(\theta,\psi) &=  11165 (2 \theta^2 + 2 \theta \psi - 7 \psi^2) 
\end{align*}
satisfying 
$f(p_1,p_2) = p_3^2$ and $g(q_1,q_2) = q_3^2$. It may then be checked
that the curve   
\[ C_3 = \left\{ 
\begin{aligned}
p_1(c,d) &= q_1(\theta,\psi) \\
p_2(c,d) &= q_2(\theta,\psi) 
\end{aligned} \right\} \subset \PP^3 \]
is everywhere locally soluble. This confirms that $\xi_1 \in S_3$.
(For the purposes of presenting this example, we adjusted the parametrisations
so that $\xi_2 = \xi_3 = 1$.)

By computing companion quadric intersections $C_2^\dagger$
and $C_3^\dagger$ as described in Section~\ref{sec:res}, we find that
$C_3$ has pushout form
\begin{align*}
 \lambda(c,d,\theta,\psi) & = 522272234 c^2 - 24659265 c d - 17398450 d^2 
  \\ & \hspace{10em} - 2641751 \theta^2 - 5524559 \theta \psi - 9670688 \psi^2. 
\end{align*}

Let $\eta \in S_3$. By the formula for the Cassels-Tate pairing in 
Section~\ref{sec:ctp}, we have $\Theta_3(\xi_1, \eta) = \sum_v 
(\lambda(P_v),\eta)_v$
where $P_v \in C_3(\Q_v)$ and $(~,~)_v$ is the 
Hilbert norm residue symbol $\Q_v^\times/(\Q_v^\times)^2 \times 
\Q_v^\times/(\Q_v^\times)^2 \to \F_2$. 
The bad primes of $E$ are 
\[ {\mathcal S} = \{ 2, 3, 5, 7, 11, 29, 71, 73, 127, 28619, 30187 \}. \]
At all primes $p \not\in {\mathcal S}$, we find that $C_3$ has good 
reduction mod $p$, whereas the two equations for $C_3$ together with 
$\lambda$ are linearly independent mod $p$. These primes therefore
make no contribution to the pairing.
Since $\eta \in (\Q^\times_v)^2$ for all $v \in \{ 71, 127, 30187, \infty \}$ 
there is also no contribution from these places.
At each of the remaining primes $p$ we find a local point
$P \in C_3(\Q_p)$ with $F(P) \equiv u_p \mod{(\Q^\times_p)^2}$ where
$u_2 = 5$ and $u_p$ is a quadratic non-residue for $p$ odd.
Therefore $\Theta_3(\xi_1, \eta) \equiv \omega(\eta) \pmod{2}$ where 
$\omega(\eta)$ is the number of prime divisors of~$\eta$. (We choose
a representative modulo squares so that $\eta$ is a square-free integer.)
This gives the first row in the following table. Repeating
for $\xi_1$ running over a basis for $S_3$ we find that
$\Theta_3 : S_3 \times S_3 \to \F_2$ is given by
\[ \begin{array}{c|ccccc}
& 15 & 73 & 87 & 231 & 28619 \\ \hline
15 & 0 & 1 & 0 & 1 & 1 \\
73 & 1 & 0 & 1 & 0 & 0 \\
87 & 0 & 1 & 0 & 0 & 1 \\
231 & 1 & 0 & 0 & 0 & 1 \\
28619 & 1 & 0 & 1 & 1 & 0
\end{array} \]
This matrix has rank $4$. Therefore $\dim S_4 = 1$
and $\rank E(\Q) = 0$. By Remark~\ref{rem:getsha} the 
$2$-primary parts of $\Sha(E/\Q)$ and 
$\Sha(E'/\Q)$ are $(\Z/4\Z)^4$ and $(\Z/2\Z)^4$.
\end{Example}

\begin{Example}
Let $E/\Q$ be the elliptic curve $y^2 = x(x^2 + a x + b)$ where
\begin{align*}
a &= -802175537664068731998722, \\
b &= 160480561352940413879437222902216664489852408321.
\end{align*}
This is an example with $E(\Q)_{\tors} \isom \Z/2\Z \times \Z/8\Z$.
Let $\phi : E \to E'$ be the $2$-isogeny with kernel generated by $T =
(0,0)$. As subgroups of $\Q^\times/(\Q^\times)^2$ we find that $S_1 =
S_2 = \langle -10, 5574 \rangle$ and
\begin{align*}
S'_1 &= \langle 3841, 920641, 262404961, 289572953761, 9289, 6049, 31441 
\rangle, \\
S'_2 &= \langle  3841, 920641, 262404961, 289572953761, 9289 \rangle. 
\end{align*}
Therefore $\rank E(\Q) \le 5$. Next we compute
$\Theta_2 : S_2 \times S'_2 \to \F_2$. The elements 
$\xi_1 = -10$ and $\eta_1 = 5574$ in $S_1 = S^{(\phi)}(E/\Q)$ 
lift to $2$-coverings of $E'$ given by
\begin{align*}
C_2 &= \{ y^2 = \lambda_1(x,z)^2 + 10 \mu_1(x,z)^2) \} \\ 
D_2 &= \{ y^2 = \lambda_2(x,z)^2 - 5574 \mu_2(x,z)^2 \}
\end{align*}
where 
\begin{align*}
\lambda_1(x,z) &= 341696479062308 (x^2 - 10 z^2) + 3516978476959251 x z, \\
\mu_1(x,z) &= 21538029761160 (x^2 - 10 z^2) + 158658854157270 x z, \\
\lambda_2(x,z) &= 1335842866662 (x^2 + 5574 z^2) + 1439937420103543 x z, \\ 
\mu_2(x,z) &= 17631567180 (x^2 + 5574 z^2) + 18855731460270 x z.
\end{align*}
It may be checked, using the formulae in \cite{testeqbq}, that $C_2$
and $D_2$ do indeed correspond to $\xi_1$ and $\eta_1$.  Since $C_2$
and $D_2$ are everywhere locally soluble, this confirms that $\xi_1$
and $\eta_1$ are in $S_2$.  Now $y - \lambda_1(x,z)$ and $y -
\lambda_2(x,z)$ are pushout forms on $C_2$ and $D_2$.  Evaluating
these at local points, and then computing sums of Hilbert norm residue
symbols, we find that $\Theta_2 : S_2 \times S'_2 \to \F_2$ is given
by
\[ \begin{array}{c|ccccc}
& 3841 & 920641 & 262404961 & 289572953761 & 9289 \\ \hline
  -10 & 0 & 0 & 0 & 0 & 1 \\
  5574 & 0 & 0 & 0 & 0 & 0 \\
\end{array} \] 
We further find that $\Theta_3$ and $\Theta'_3$ are identically zero.
Therefore $S_3 = S_4 = \langle 5574 \rangle$ and 
\[ S'_3 = S'_4 = \langle 3841, 920641, 262404961, 289572953761
\rangle. \]

Using the methods in Sections~\ref{sec:pushout} and~\ref{sec:res} we
compute an everywhere locally soluble $2$-covering $D_4$ of $D_2$ (and
hence $4$-covering of $E'$) with equations
\begin{align*}
& 9055 x_1^2 + 139619 x_1 x_2 + 394387 x_1 x_3
 + 47027 x_1 x_4 - 94269 x_2^2  \\
& \hspace{1em} - 234422 x_2 x_3 + 266438 x_2 x_4 + 127750 x_3^2
 - 130775 x_3 x_4 - 137150 x_4^2
 = 0, \\ 
& 255171 x_1^2 - 961185 x_1 x_2 + 297383 x_1 x_3
 + 224191 x_1 x_4 + 152028 x_2^2
  \\
& \hspace{1em} + 461745 x_2 x_3 - 59967 x_2 x_4 + 370158 x_3^2
 - 350350 x_3 x_4 - 198938 x_4^2 
=0.
\end{align*}
We find that $D_4$ has pushout form
\begin{align*}
\lambda &= 2471492850764286 x_1^2 - 573148078730175 x_1 x_2
 + 136617115364660 x_1 x_3 \\
& - 1043769539460119 x_1 x_4 - 338965008210675 x_2^2
 + 46036721190885 x_2 x_3 \\
& + 1057225484721135 x_2 x_4 - 989678424716819 x_3^2
 + 2171876872481818 x_3 x_4 \\ 
& - 1169960121148148 x_4^2.
\end{align*}
Using this we compute that 
$\Theta_4 : S_4 \times S'_4 \to \F_2$ is given by
\[ \begin{array}{c|ccccc}
&  3841 & 920641 & 262404961 & 289572953761 \\ \hline
  5574 & 1 & 1 & 0 & 1 
\end{array} \] 
Therefore $S_5 = 0$, $\dim S'_5 = 3$ and
$\rank E(\Q) \le 1$. If $\rank E(\Q) = 1$ then by
Remark~\ref{rem:getsha} the $2$-primary parts of $\Sha(E/\Q)$ and
$\Sha(E'/\Q)$ are $(\Z/2\Z)^2 \times (\Z/4\Z)^2$ and
 $(\Z/2\Z)^4 \times (\Z/4\Z)^2$.
\end{Example}

\begin{Example}
Let $E/\Q$ be the elliptic curve $y^2 = x^3 - d^2 x$ where
\[d = 743114132612994 = 2 \times 3 \times 19 \times 953 \times 1427
\times 2137 \times 2243.\] This was one of the candidates in
\cite{Watkins+} for a congruent number elliptic curve of rank $6$.
Let $\phi : E \to E'$ be the $2$-isogeny with kernel generated by $T =
(-d,0)$. We find that
\begin{align*}
  S_1 &= S_2 = S_3 = S_4 =\langle 1906, 2137 \rangle \subset \Q^\times/(\Q^\times)^2 , \\
  S'_1 &= S'_2 = S'_3 = S'_4 = \langle 2, 57, 953, 2137, 4281, 6729
  \rangle \subset \Q^\times/(\Q^\times)^2 .
\end{align*}
Therefore $\rank E(\Q) \le 6$. (We note that $1906 = 2 \times 953$,
$57 = 3 \times 19$, $4281 = 3 \times 1427$ and $6729 = 3 \times 2243$.)
The elements $1906$
and $2137$ in $S^{(\phi)}(E/\Q)$
lift to $4$-coverings of $E'$ with equations 
\begin{align*}
&  1188 x_1^2 + 1244 x_1 x_2 + 732 x_1 x_3 + 5599 x_1 x_4
   - 1530 x_2^2  \\ & \hspace{1em} - 3687 x_2 x_3 + 2824 x_2 x_4 - 2928 x_3^2 
    + 1780 x_3 x_4 + 4886 x_4^2 = 0, \\
&  3298 x_1^2 - 7382 x_1 x_2 + 3881 x_1 x_3 - 3470 x_1 x_4 
   + 20 x_2^2  \\ & \hspace{1em} + 2136 x_2 x_3 + 7147 x_2 x_4 
   + 1517 x_3^2 - 4464 x_3 x_4 - 2455 x_4^2  = 0, 
\end{align*}
and 
\begin{align*}
& 68 x_1^2 + 416 x_1 x_2 + 422 x_1 x_3 + 1372 x_1 x_4 
  + 401 x_2^2  \\ & \hspace{1em} + 1146 x_2 x_3 + 1154 x_2 x_4 + 475 x_3^2 
  + 1528 x_3 x_4 - 2366 x_4^2 = 0, \\
& 18869 x_1^2 - 13870 x_1 x_2 - 14599 x_1 x_3 - 2402 x_1 x_4 
  + 2322 x_2^2  \\ & \hspace{1em} + 15142 x_2 x_3 - 1629 x_2 x_4 - 12250 x_3^2 
   + 6366 x_3 x_4 + 4913 x_4^2 = 0.
\end{align*}
These have pushout forms
\begin{align*}
\lambda &= 280117153304 x_1^2 + 627376555572 x_1 x_2 - 534852420548 x_1 x_3 \\ & - 
        24376482389 x_1 x_4 + 7830950834 x_2^2 - 165910883299 x_2 x_3 \\ & + 
        255303594426 x_2 x_4 - 574956757207 x_3^2 - 923597568302 x_3 x_4 \\ & - 
        398161865115 x_4^2,
\end{align*}
and
\begin{align*}
\lambda &= 302992919 x_1^2 + 165225436 x_1 x_2 - 554084592 x_1 x_3 
   - 3000363446 x_1 x_4 \\ & + 259816709 x_2^2 - 1891629098 x_2 x_3 
    + 2507867060 x_2 x_4 - 1337797147 x_3^2 \\ & + 572269312 x_3 x_4 
    - 1348920925 x_4^2.
\end{align*}

Using these we compute that $\Theta_4 : S_4 \times S'_4 \to \F_2$
is given by
\[ \begin{array}{c|cccccc}
& 2 & 57 & 953 & 2137 & 4281 & 6729 \\ \hline
1906 & 0 & 1 & 1 & 0 & 0 & 0 \\
2137 & 0 & 0 & 1 & 0 & 1 & 0
\end{array} \]
This matrix has rank $2$. Therefore
$S_5=0$, $\dim S'_5 = 4$ and
$\rank E(\Q) \le 2$. Searching for rational points on 
$2$-coverings of $E$, and on $2$-coverings of elliptic curves isogenous to $E$,
we find that $E(\Q)$ has independent points of infinite order
\begin{align*}
  P_1 & = (378998487378128086417/42^2,7378254260469086949865651299481/42^3), \\
P_2 & = (4906543485739785817699690091987274941361019066993/76947716\backslash \\
& \hspace{3em} 
   400301612^2, 480965296472375270212047079501479961953600924\backslash \\
& \hspace{3em} 
  7328327459226000929166668023/76947716400301612^3).
\end{align*}
Therefore $\rank E(\Q) = 2$.
By Remark~\ref{rem:getsha} the
$2$-primary parts of $\Sha(E/\Q)$ and $\Sha(E'/\Q)$ 
are isomorphic to $(\Z/4\Z)^4$. Similar 
calculations show that the other two elliptic 
curves in the isogeny class have $2$-primary 
parts of $\Sha$ isomorphic to 
$(\Z/4\Z)^2 \times (\Z/8\Z)^2$ and $(\Z/8\Z)^4$.
\end{Example}


\end{document}